\documentclass[a4paper,leqno,twoside]{article}
\usepackage{latexsym,amssymb,amsmath, amsthm}
\usepackage[affil-it]{authblk}

\newtheorem{theo}{Theorem}[section]
\newtheorem{coro}[theo]{Corollary}
\newtheorem{lemm}[theo]{Lemma}
\newtheorem{prop}[theo]{Proposition}

\newtheorem{exam} [theo]{Example}
\newtheorem{rema}[theo]{Remark}

\author{Ivan D. Chipchakov\\
Institute of Mathematics and Informatics, Bulgarian Academy\\ of 
Sciences, Acad. G. Bonchev Str., bl. 8, 1113, Sofia, Bulgaria}
\title{On index-exponent relations over Henselian
fields with local residue fields}
\date{}
\begin{document}
\maketitle

\begin{abstract}
Let $p$ be a prime and $(K, v)$ a Henselian valued field with
a residue field $\widehat K$. This paper determines the Brauer
$p$-dimension of $K$, in case $p \neq {\rm char}(\widehat K)$ and
$\widehat K$ is a $p$-quasilocal field properly included in its
maximal $p$-extension. When $\widehat K$ is a local field, it
describes index-exponent pairs of central division $K$-algebras of
$p$-primary degrees. The same goal is achieved, if $(K, v)$ is
maximally complete, char$(K) = p$ and $\widehat K$ is local.
\end{abstract}


\par
\section{\bf Introduction}
\par
\medskip
Let $E$ be a field, $\mathbb P$ the set of prime numbers, and for
each $p \in \mathbb P$, let $E(p)$ be the maximal $p$-extension of
$E$ in a separable closure $E _{\rm sep}$, and $r _{p}(E)$ the rank
of the Galois group $\mathcal{G}(E(p)/E)$ as a pro-$p$-group (put $r
_{p}(E) = 0$, if $E(p) = E$). Denote by $s(E)$ the class of
finite-dimensional associative central simple $E$-algebras, and by
$d(E)$ the subclass of division algebras $D \in s(E)$. For each $A
\in s(E)$, let $[A]$ be the equivalence class of $A$ in the Brauer
group Br$(E)$, and $D _{A}$ a representative of $[A]$ lying in
$d(E)$. The existence of $D _{A}$ and its uniqueness, up-to an
$E$-isomorphism, is established by Wedderburn's structure theorem
(cf. \cite{P}, Sect. 3.5), which implies the dimension $[A\colon E]$
is a square of a positive integer deg$(A)$ (the degree of $A$). It
is known that Br$(E)$ is an abelian torsion group, so it decomposes
into the direct sum, taken over $\mathbb P$, of its $p$-components
Br$(E) _{p}$ (see \cite{P}, Sects. 3.5 and 14.4). The Schur index
ind$(D) = {\rm deg}(D _{A})$ and the exponent exp$(A)$, i.e. the
order of $[A]$ in Br$(E)$, are invariants of both $A$ and $[A]$.
Their general relations and behaviour under scalar extensions of
finite degrees are described as follows (cf. \cite{P}, Sects. 13.4,
14.4 and 15.2):
\par
\medskip
(1.1) (a) exp$(A) \mid {\rm ind}(A)$ and $p \mid {\rm exp}(A)$, for
each $p \in \mathbb P$ dividing ind$(A)$. For any $B \in s(E)$ with
ind$(B)$ prime to ind$(A)$, ind$(A \otimes _{E} B) = {\rm
ind}(A).{\rm ind}(B)$; if $A$, $B \in d(E)$, then the tensor product
$A \otimes _{E} B$ lies in $d(E)$;
\par
(b) ind$(A)$ and ind$(A \otimes _{E} R)$ divide ind$(A \otimes _{E}
R)[R\colon E]$ and ind$(A)$, respectively, for each finite field
extension $R/E$ of degree $[R\colon E]$.
\par
\medskip
As shown by Brauer (see, e.g., \cite{P}, Sect. 19.6), (1.1) (a)
determines all generally valid index-exponent relations. It is
known, however, that, for a number of fields $E$, the pairs ind$(A),
{\rm exp}(A)$, $A \in s(E)$, are subject to much tougher
restrictions than those described by (1.1) (a). The Brauer
$p$-dimensions Brd$_{p}(E)$, $p \in \mathbb P$, contain essential
information on these restrictions. We say that Brd$_{p}(E) = n$,
where $n \in \mathbb Z$, if $n$ is the least integer $\ge 0$ for
which ind$(D) \le {\rm exp}(D) ^{n}$ whenever $D \in d(E)$ and
$[D] \in {\rm Br}(E) _{p}$; if no such $n$ exists, we put Brd$_{p}(E)
= \infty $. In view of (1.1), Brd$_{p}(E) \le 1$, for a given $p$, if
and only if ind$(D) = {\rm exp}(D)$, for each $D \in d(E)$ with $[D]
\in {\rm Br}(E) _{p}$; Brd$_{p}(E) = 0$ if and only if Br$(E) _{p} =
\{0\}$. The absolute Brauer $p$-dimension abrd$_{p}(E)$ of $E$ is
defined as the supremum Brd$_{p}(R)\colon R \in {\rm Fe}(E)$,
Fe$(E)$ being the set of finite extensions of $E$ in $E _{\rm sep}$.
For example, when $E$ is a global or local field, Brd$_{p}(E) = {\rm
abrd}_{p}(E) = 1$, $p \in \mathbb P$, and there exist $Y _{n} \in
d(E)$, $n \in \mathbb N$, with ind$(Y _{n}) = n$, for any $n$ (see
\cite{We}, Ch. XII, Sect. 2; Ch. XIII, Sects. 3, 6).
\par
\medskip
This paper deals with the study of index-exponent $K$-pairs, for a
Henselian (valued) field $(K, v)$, along the lines drawn in
\cite{Ch5}, Sect. 5. Its purpose is to determine Brd$_{p}(K)$ and to
describe $p$-primary index-exponent $K$-pairs, provided that
the residue field $\widehat K$ of $(K, v)$ is endowed with a
Henselian discrete valuation $\omega $ whose residue field  is
quasifinite, and $p \in \mathbb P$ is different from char$(\widehat
K)$ (for other types of $\widehat K$, such as the one of a global
field, see \cite{Ch5}, Sect. 5). Our main result, presented by the
following theorem, concerns the case where $\widehat K$ is a local
field and the value group $v(K)$ is $p$-indivisible, i.e. the
quotient group $v(K)/pv(K)$ is nontrivial. When $K$ contains a
primitive $p$-th root of unity, it shows that index-exponent
$p$-primary $K$-pairs are not determined only by Brd$_{p}(K)$:

\medskip
\begin{theo}
\label{theo1.1} Let $(K, v)$ be a Henselian field with {\rm
Brd}$_{p}(K) < \infty $, for some
\par\noindent
$p \in \mathbb P$, and let $m _{p} = {\rm min}\{\tau (p), r
_{p}(\widehat K)\}$, where $\tau (p)$ is the dimension of
\par\noindent $v(K)/pv(K)$ as a vector space over the field
$\mathbb F _{p} = \mathbb Z/p\mathbb Z$. Assume that $\tau (p) > 0$,
$p \neq {\rm char}(\widehat K)$, $\widehat K$ is a local field, and
$\varepsilon _{p}$ is a primitive $p$-th root of unity in $\widehat K
_{\rm sep}$, denote by $\nu $ the greatest integer for which
$\widehat K$ contains a primitive $p ^{\nu }$-th root of unity, and
in case $\varepsilon _{p} \in \widehat K$, put $r _{p} ^{\prime
}(\widehat K) = r _{p}(\widehat K) - 1$ and $m _{p} ^{\prime } = {\rm
min}\{\tau (p), r _{p} ^{\prime }(\widehat K)\}$. For each $n \in
\mathbb N$, let $\mu (p, n) = nm _{p} ^{\prime } + \nu _{n}(m _{p} -
m _{p} ^{\prime } + [(\tau (p) - m _{p})/2])$, if $\varepsilon _{p}
\in \widehat K$, where $\nu _{n} = {\rm min}\{n, \nu \}$, and $\mu
(p, n) = nm _{p}$, if $\varepsilon _{p} \notin \widehat K$. Then {\rm
Brd}$_{p}(K) = \mu (p, 1)$; moreover, for a pair $(k, n) \in \mathbb N
^{2}$, there exists $D _{k,n} \in d(K)$ with {\rm ind}$(D _{k,n}) =
p ^{k}$ and {\rm exp}$(D _{k, n}) = p ^{n}$ if and only if $n \le k
\le \mu (p, n)$.
\end{theo}
\par
\medskip
Assuming that $(K, v)$ is Henselian, $p \in \mathbb P$ is not equal
to char$(\widehat K)$, $\tau (p)$ and $\varepsilon _{p}$ are
defined as above, and $(\widehat K, \omega )$ is a Henselian
discrete valued field (abbr, an HDV-field) with a quasifinite residue field
$\tilde k$, we obtain the following result:
\par
\medskip
(1.2) (a) If $0 < \tau (p) < \infty $, char$(\widehat K) = 0$, and
$\tilde k$ is infinite with char$(\tilde k) = p$, then Brd$_{p}(K) =
\tau (p)$ and $(p ^{k}, p ^{n})$, $k, n \in \mathbb N$, $n \le k \le
n\tau (p)$, are all nontrivial index-exponent $p$-primary $K$-pairs;
\par
(b) Brd$_{p}(K) = 1$ and $(p ^{n}, p ^{n})$, $n \in \mathbb N \cup
\{0\}$, are all index-exponent $K$-pairs, in case $p \neq {\rm
char}(\tilde k)$ and $\varepsilon _{p} \notin \widehat K$; the same
holds, if $p \neq {\rm char}(\tilde k)$ and $\tau (p) \le 1$;
\par
(c) When $p \neq {\rm char}(\tilde k)$, $\varepsilon _{p} \in
\widehat K$, and $2 \le \tau (p) < \infty $, we have $r _{p}(\widehat 
K) = 2$ and Brd$_{p}(K) = 1 + [\tau (p)/2]$;
\par
(d) In the setting of (c), if $\widehat K$ contains finitely many 
roots of unity of $p$-primary degrees, then index-exponent 
$p$-primary $K$-pairs are determined in accordance with Theorem 
\ref{theo1.1}; when $\widehat K$ contains infinitely many such roots, 
$(p ^{k}, p ^{n})$, $k, n \in \mathbb N$, $n \le k \le n{\rm 
Brd}_{p}(K)$, are index-exponent $K$-pairs.
\par
\medskip
When $(K, v)$ is a maximally complete field with char$(K) = p$ and
$\widehat K$ a local field, Brd$_{p}(K)$ and index-exponent
$p$-primary $K$-pairs are determined as follows:
\par
\medskip
\begin{prop}
\label{prop1.2} Assume that $(K, v)$ is a maximally complete field,
{\rm char}$(K) = p > 0$, and $\widehat K$ is a local field, and 
define $\tau (p)$ as in Theorem \ref{theo1.1}. Then:
\par
{\rm (a)} {\rm Brd}$_{p}(K) = \infty $ if and only if $\tau (p) =
\infty $; when this holds, $(p ^{k}, p ^{n})$ is an index-exponent
pair over $K$, for any $k, n \in \mathbb N$ with $k \ge n$;
\par
{\rm (b)} {\rm Brd}$_{p}(K) = \tau (p)$, provided that $\tau (p) <
\infty $; in this case, $(p ^{k}, p ^{n})$ is an index-exponent
$K$-pair, where $k, n \in \mathbb N$, if and only if $n \le k \le
n\tau (p)$.
\end{prop}

\medskip
Proposition \ref{prop1.2} is deduced in Section 3 from our description 
of index-exponent $p$-primary pairs over maximally complete fields of 
characteristic $p$ with perfect residue fields (see Corollary 
\ref{coro3.6} and Proposition \ref{prop3.5}). The proofs of (1.2) and 
Theorem \ref{theo1.1} rely on the fact that HDV-fields with 
quasifinite residue fields are quasilocal, i.e. their finite 
extensions are $p$-quasilocal fields with respect to every $p \in 
\mathbb P$ (see \cite{S1}, Ch. XIII, Sect. 3). As in \cite{Ch2}, a 
field $E$ with $r _{p}(E) > 0$, for some $p$, is called 
$p$-quasilocal, if the relative Brauer group Br$(E ^{\prime }/E)$ 
equals the group $_{p}{\rm Br}(E) = \{b \in {\rm Br}(E)\colon pb = 
0\}$, for every degree $p$ extension $E ^{\prime }$ of $E$ in $E(p)$; 
when $r _{p}(E) = 0$, we say that $E$ is $p$-quasilocal, if Br$(E) 
_{p} = \{0\}$. The part of Theorem \ref{theo1.1} concerning 
Brd$_{p}(K)$ is a special case of a formula for Brd$_{p}(K)$, deduced 
when $\widehat K$ is any $p$-quasilocal field with char$(\widehat K) 
\neq p$ and $r _{p}(\widehat K) > 0$ (see Section 4). To prove this 
formula we use the inequality Brd$_{p}(\widehat K) \le 1$, the 
surjectivity of the scalar extension map Br$(\widehat K) _{p} \to 
{\rm Br}(\widehat K ^{\prime }) _{p}$, for every extension $\widehat 
K ^{\prime }$ of $\widehat K$ in $\widehat K(p)$, and the following 
relations between finite extensions of $\widehat K$ in $\widehat 
K(p)$ and algebras $\Delta _{p} \in d(\widehat K)$ of $p$-primary 
degrees (see \cite{Ch2}, I, Sects. 3 and 4):
\par
\medskip
(1.3) ind$(\Delta _{p}) = {\rm g.c.d.}\{[L _{p}\colon \widehat K],
{\rm ind}(\Delta _{p})\}{\rm ind}(\Delta _{p} \otimes _{\widehat K}
L _{p})$ whenever $L _{p}$ is a finite extension of $\widehat K$ in
$\widehat K(p)$. Specifically, $L _{p}$ embeds in $\Delta _{p}$ as a
$\widehat K$-subalgebra if and only if $[L _{p}\colon \widehat K]
\mid {\rm ind}(\Delta _{p})$; $[\Delta _{p}] \in {\rm Br}(L
_{p}/\widehat K)$ if and only if ind$(\Delta _{p}) \mid [L
_{p}\colon \widehat K]$.
\par
\medskip
Statements (1.2) and the concluding assertion of Theorem 
\ref{theo1.1} are proved in Section 5. Their proofs are based on 
Morandi's theorem \cite{Mo}, the theory of division algebras over 
Henselian fields developed in \cite{JW}, and the structure of the 
(continuous) character group $C(\widehat K(p)/\widehat K)$ of
$\mathcal{G}(\widehat K(p)/\widehat K)$ as an abstract abelian group
(see (5.2), (5.3) and Remark \ref{rema5.3}). Our proofs also rely on
the fact that if $\widehat K$ is a local field or $p \neq {\rm char}(\tilde
k)$, then $\mathcal{G}(\widehat K(p)/\widehat K)$ is a Demushkin
group if $\widehat K$ contains a primitive $p$-th root of unity, and
$\mathcal{G}(\widehat K(p)/\widehat K)$ is a free pro-$p$-group,
otherwise (cf. \cite{S2}, Ch. II, 2.2 and 5.6). By a Demushkin group,
we mean a pro-$p$-group $G _{p}$ whose continuous cohomology
groups $H ^{i}(G _{p}, \mathbb F _{p})$ with coefficients in $\mathbb
F _{p}$, for $i = 1, 2$, satisfy the following conditions: $H ^{2}(G
_{p}, \mathbb F _{p})$ is of order $p$, $H ^{1}(G _{p}, \mathbb F
_{p})$ is finite and abelian of period $p$, and for any nonzero $a
\in H ^{1}(G _{p}, \mathbb F _{p})$, the homomorphism $\varphi
_{a}\colon H ^{1}(G _{p}, \mathbb F _{p}) \to H ^{2}(G _{p}, \mathbb
F _{p})$, mapping each $b \in H ^{1}(G _{p}, \mathbb F _{p})$ into
the cup-product $a \cup b$, is surjective. We also use the well-known 
fact that local fields contain finitely many roots of unity, and take 
into account that Brauer groups of HDV-fields with quasifinite 
residue fields are isomorphic to the quotient group $\mathbb 
Q/\mathbb Z$ of the additive group of rational numbers by the 
subgroup of integers (cf. \cite{S1}, Ch. XIII, Sect. 3).

\medskip
The basic notation and terminology used and conventions kept in this
paper are standard, like in \cite{Ch2}, I, \cite{Ch4} and \cite{Ch5}. 
We write $Z(B)$ for the centre of an associative ring $B$. Given a 
Henselian field $(K, v)$, $K _{\rm ur}$ denotes the compositum of 
inertial extensions of $K$ in $K _{\rm sep}$; the notions of an 
inertial, a nicely semi-ramified (abbr, NSR), and a totally ramified 
(division) $K$-algebra, are defined in \cite{JW}. Section 2 includes
valuation-theoretic preliminaries used in the sequel. By a
Pythagorean field, we mean a formally real field whose set of squares
is additively closed. As usual, $[r]$ stands for the integral part of a
real number $r \ge 0$, and for any $p \in \mathbb P$, a $\mathbb Z 
_{p}$-extension means a Galois extension whose Galois group is 
isomorphic to the additive group $\mathbb Z _{p}$ of $p$-adic 
integers. The set of intermediate fields of a field extension 
$\Lambda /\Psi $ is denoted by $I(\Lambda /\Psi )$. Symbol algebras 
are defined, e.g., in \cite{JW} and \cite{P}, Sect. 15.4. Galois 
groups are viewed as profinite with respect to the Krull topology, 
and by a profinite group homomorphism, we mean a continuous one. The 
reader is referred to \cite{L}, \cite{Efr2}, \cite{JW}, \cite{FV}, 
\cite{P} and \cite{S2}, for missing definitions concerning field 
extensions, orderings and valuations, $m$-dimensional local fields, 
simple algebras, Brauer groups and Galois cohomology.

\medskip
\section{\bf Preliminaries}
\par
\medskip
Let $(K, v)$ be a Krull valued field with a residue field $\widehat
K$ and a (totally ordered) value group $v(K)$. We say that $(K, v)$
is Henselian, if $v$ extends uniquely, up-to an equivalence, to a
valuation $v _{L}$ on each algebraic extension $L/K$. This occurs,
for example, if $(K, v)$ is maximally complete, i.e. it has no 
immediate proper extension (a valued extension $(K ^{\prime }, v 
^{\prime })$, such that $K ^{\prime } \neq K$, $\widehat K ^{\prime } 
= \widehat K$ and $v ^{\prime }(K ^{\prime }) = v(K)$). When $(K, v)$ 
is Henselian, we denote by $\widehat L$ the residue field of $(L, v 
_{L})$ and put $v(L) = v _{L}(L)$, for any algebraic extension $L/K$. 
Clearly, $\widehat L/\widehat K$ is an algebraic extension and $v(K)$ 
is an ordered subgroup of $v(L)$; $e(L/K)$ denotes the index of 
$v(K)$ in $v(L)$. By Ostrowski's theorem (cf. \cite{Efr2}, 
Theorem~17.2.1), when $L/K$ is finite, $[L\colon K]$, $[\widehat 
L\colon \widehat K]$ and $e(L/K)$ are related as follows:
\par
\medskip
(2.1) $[\widehat L\colon \widehat K]e(L/K)$ divides $[L\colon K]$ and
$[L\colon K][\widehat L\colon \widehat K] ^{-1}e(L/K) ^{-1}$ is not
divisible by any $p \in \mathbb P$, $p \neq {\rm char}(\widehat K)$;
$[L\colon K] = [\widehat L\colon \widehat K]e(L/K)$, if
char$(\widehat K) \nmid [L\colon K]$.
\par
\medskip
The Henselity of $(K, v)$ ensures that each $\Delta \in d(K)$ has a
unique, up-to an equivalence, valuation $v _{\Delta }$ extending $v$
and possessing an abelian value group $v(\Delta )$ (cf. \cite{Sch},
Ch. 2, Sect. 7). This group is totally ordered and includes $v(K)$ as
an ordered subgroup of index $e(\Delta /K) \le [\Delta \colon K]$.
Also, the residue division ring $\widehat \Delta $ of $(\Delta , v
_{\Delta })$ is a $\widehat K$-algebra, and by Ostrowski-Draxl's
theorem \cite{Dr2}, $e(\Delta /K)[\widehat \Delta \colon \widehat K]
\mid [\Delta \colon K]$ and if char$(\widehat K)$ $\nmid $
ind$(\Delta )$, then $[\Delta \colon K] = e(\Delta /K)[\widehat
\Delta \colon \widehat K]$. Statement (2.1) and the Henselity of $(K,
v)$ imply the following:
\par
\medskip
(2.2) The quotient groups $v(K)/pv(K)$ and $v(L)/pv(L)$ are
isomorphic, if $p \in \mathbb P$ and $[L\colon K] < \infty $. When
char$(\widehat K) \nmid [L\colon K]$, the natural embedding of $K$
into $L$ induces canonically an isomorphism $v(K)/pv(K) \cong
v(L)/pv(L)$.
\par
\medskip
A finite extension $R$ of $K$ is said to be inertial, if $[R\colon
K] = [\widehat R\colon \widehat K]$ and $\widehat R/\widehat K$ is
separable. We say that $R/K$ is totally ramified, if $[R\colon K] =
e(R/K)$; $R/K$ is called tamely ramified, if $\widehat
R/\widehat K$ is separable and char$(\widehat K) \nmid e(R/K)$. The
properties of $K _{\rm ur}/K$ used in the sequel are essentially those
presented in \cite{JW}, page 135, and restated in \cite{Ch3}, (3.3) 
(see also \cite{TW}, Theorem A.24). Here we recall some results on 
central division $K$-algebras (most of which can be found in 
\cite{JW}):
\par
\medskip
(2.3) (a) If $D \in d(K)$ and char$(\widehat K) \nmid {\rm
ind}(D)$, then $[D] = [S \otimes _{K} V \otimes _{K} T]$, for some
$S$, $V$, $T \in d(K)$, such that $S/K$ is inertial, $V/K$ is NSR,
$T/K$ is totally ramified, $T \otimes _{K} K _{\rm ur} \in d(K _{\rm
ur})$, exp$(T \otimes _{K} K _{\rm ur}) = {\rm exp}(T)$, and $T$ is
a tensor product of totally ramified cyclic $K$-algebras (see also
\cite{Dr2}, Theorem~1);
\par
(b) The set IBr$(K) = \{[S ^{\prime }] \in {\rm Br}(K)\colon S
^{\prime } \in d(K), S ^{\prime }/K$ inertial$\}$ is a subgroup of
Br$(K)$ canonically isomorphic to Br$(\widehat K)$; Brd$_{p}(\widehat
K) \le {\rm Brd}_{p}(K)$, $p \in \mathbb P$, and equality holds when
$p \neq {\rm char}(\widehat K)$ and $v(K) = pv(K)$;
\par
(c) With assumptions and notation being as in (a), if $T \neq K$,
then $K$ contains a primitive root of unity of degree exp$(T)$; in
addition, if $T _{n} \in d(K)$ and $[T _{n}] = n[T] \neq 0$, for some
$n \in \mathbb N$, then $T _{n}/K$ is totally ramified;

\medskip
Statement (2.3) can be supplemented as follows (see, e.g.,
\cite{Ch5}, Sect. 4):
\par
\medskip
(2.4) If $D$, $S$, $V$ and $T$ are related as in (2.3) (a), then:
\par
(a) $n[D] \in {\rm IBr}(K)$, for a given $n \in \mathbb N$, if and
only if exp$(V) \mid n$ and exp$(T) \mid n$;
\par
(b) $D/K$ is inertial if and only if $V = T = K$; $D/K$ is
inertially split, i.e. $[D] \in {\rm Br}(K _{\rm ur}/K)$, if and
only if $T = K$;
\par
(c) exp$(D) = {\rm lcm}({\rm exp}(S), {\rm exp}(V), {\rm exp}(T))$.

\par
\medskip\noindent
The following result of \cite{Ch5} gives a formula for Brd$_{p}(K)$
whenever $p \neq {\rm char}(\widehat K)$ and Brd$_{p}(\widehat K) =
0$:

\medskip
\begin{theo}
\label{theo2.1} Assume that $(K, v)$ is a Henselian field with {\rm
Brd}$_{p}(\widehat K) < \infty $, for some $p \in \mathbb P$, $p
\neq {\rm char}(\widehat K)$, and let $\tau (p)$, $\varepsilon
_{p}$ and $m _{p}$ be as in Theorem \ref{theo1.1}. Then:
\par
{\rm (a)} {\rm Brd}$_{p}(K) = \infty $ if and only if $m _{p} =
\infty $ or $\tau (p) = \infty $ and $\varepsilon _{p} \in \widehat
K$;
\par
{\rm (b)} $[(\tau (p) + m _{p})/2] \le {\rm Brd}_{p}(K) \le {\rm
Brd}_{p}(\widehat K) + [(\tau (p) + m _{p})/2]$, if $\tau
(p) < \infty $ and $\varepsilon _{p} \in \widehat K$; when $m _{p} <
\infty $ and $\varepsilon _{p} \notin \widehat K$, $m _{p} \le {\rm
Brd}_{p}(K) \le {\rm Brd}_{p}(\widehat K) + m _{p}$.
\end{theo}

\medskip
As shown in \cite{Ch5}, Sect. 4, Theorem \ref{theo2.1} leads to the
following description of index-exponent $p$-primary $K$-pairs, in
the case where Brd$_{p}(K) = \infty $:

\medskip
\begin{coro}
\label{coro2.2} Let $(K, v)$ be a Henselian field with {\rm
Brd}$_{p}(\widehat K) < \infty  = {\rm Brd}_{p}(K)$, for some $p \neq
{\rm char}(\widehat K)$. Then the following alternative holds:
\par
{\rm (a)} $(p ^{k}, p ^{n})\colon k, n \in \mathbb N, n \le k$, are
index-exponent $K$-pairs;
\par
{\rm (b)} $p = 2$ and $\widehat K$ is a Pythagorean field; such
being the case, the group {\rm Br}$(K) _{2}$ has period $2$, and
there are $D _{m} \in d(K)$, $m \in \mathbb N$, with {\rm ind}$(D
_{m}) = 2 ^{m}$.
\end{coro}

\medskip
This Section ends with a lemma that is implicitly used in the proofs
of the main results of the following Section.

\medskip
\begin{lemm}
\label{lemm2.3} Let $(K, v)$ be a valued field with {\rm char}$(K) =
p > 0$ and $v(K) \neq pv(K)$, and let $\pi \in K ^{\ast }$ be an
element of value $v(\pi ) \notin pv(K)$. Assume that $G$ is a finite
$p$-group of order $p ^{t}$. Then there exists a Galois extension $M$
of $K$ in $K(p)$, such that $\mathcal{G}(M/K) \cong G$, $v$ is
uniquely, up-to an equivalence, extendable to a valuation $v _{M}$ of
$M$, and $v(\pi ) \in p ^{t}v _{M}(M)$; in particular, $v
_{M}(M)/v(K)$ is cyclic and $(M, v _{M})/(K, v)$ is totally ramified.
\end{lemm}

\medskip
\begin{proof}
One may assume, for the proof, that $v(\pi ) < 0$. Let $\mathbb F$ be
the prime subfield of $K$, $(K _{v}, \bar v)$ a Henselization of $(K, v)$,
$\rho (K _{v}) = \{u ^{p} - u\colon \ u \in K _{v}\}$, $\omega $
the valuation of the field $\Phi = \mathbb F(\pi ) $ induced by $v$
and for each $m \in \mathbb N$, let $L _{m}$ and $\Lambda _{m}$ be the
root fields in $K _{\rm sep}$ over $K$ and $\Phi $, respectively, of
the polynomial $f _{m}(X) = X ^{p} - X - \pi _{m}$, where $\pi _{m} =
\pi ^{1+qm}$. Identifying $K _{v}$ with its $K$-isomorphic copy in $K
_{\rm sep}$, take a Henselization $(\Phi _{\omega }, \bar \omega )$
of $(\Phi , \omega )$ among the valued subfields of $(K _{v}, \bar v)$
(this is possible, by \cite{Efr2}, Theorem~15.3.5), and put
\par\noindent $\Psi _{m}
= \Lambda _{1} \dots \Lambda _{m}$ and $M _{m} = L _{1} \dots L
_{m}$, for each $m$. It is well-known that $(K _{v}, \bar v)/(K, v)$ and
$(\Phi _{\omega }, \bar \omega )/(\Phi , \omega )$ are immediate
extensions, i.e. 
$\widehat K _{v} = \widehat K$, $\bar v(K _{v}) = v(K)$ and $\widehat
\Phi _{\omega } = \widehat \Phi $, $\bar \omega (\Phi _{\omega }) =
\omega (\Phi )$. Also, it is easily verified that $\rho (K _{v})$ is
an $\mathbb F$-subspace of $K _{v}$, and $\bar v(u ^{\prime }) \in
pv(K)$ whenever $u ^{\prime } \in \rho (K _{v})$ and $\bar v(u
^{\prime }) < 0$. This implies the cosets $\pi _{m} + \rho (K _{v})$,
$m \in \mathbb N$, are linearly independent over $\mathbb F$, so the
Artin-Schreier theorem (cf. \cite{L}, Ch. VIII, Sect. 6) enables one
to prove the following statement, for each $m \in \mathbb N$:
\par
\medskip
(2.5) $L _{m}/K$, $L _{m}K _{v}/K _{v}$, $\Lambda _{m}/\Phi $ and
$\Lambda _{m}\Phi _{\omega }/\Phi _{\omega }$ are degree $p$ cyclic
extensions; $M _{m}/K$, $M _{m}K _{v}/K _{v}$, $\Psi _{m}/\Phi $ and
$\Psi _{m}\Phi _{\omega }/\Phi _{\omega }$ are abelian of degree $p
^{m}$.
\par
\medskip\noindent
Let now $G _{r}$ be a finite $p$-group of rank $r > 0$ and order $p
^{\mu (r)}$. Since char$(\Phi ) = p$, and therefore,
$\mathcal{G}(\Phi (p)/\Phi )$ is a free pro-$p$-group (cf. \cite{S2},
I, 1.5, 4.2; II, 2.2), there exists a Galois extension $\Gamma _{r}$
of $\Phi $ in $K _{\rm sep}$, such that $\mathcal{G}(\Gamma _{r}/\Phi
) \cong G _{r}$ and $\Psi _{r} \in I(\Gamma _{r}/\Phi )$. Hence, by
Galois theory, the field $\Gamma _{r}K$ is a Galois extension of $K$
with $\mathcal{G}(\Gamma _{r}K/K) \cong \mathcal{G}(\Gamma _{r}/\Phi
) \cong G _{r}$. We prove that $\Gamma _{r}K/K$, $G _{r}$ and $\pi $
are related in agreement with Lemma \ref{lemm2.3}. Firstly, it is
easy to see that $\Psi _{r}$ equals the fixed field of the Frattini
subgroup of $\mathcal{G}(\Gamma _{r}/\Phi )$. Secondly, it follows
from the Artin-Schreier theorem and the definition of $\Psi _{r}$
that every degree $p$ extension of $\Phi _{\omega }$ in $\Psi
_{r}\Phi _{\omega }$ is totally ramified (relative to $\bar \omega
$). Note also that $\widehat \Phi $ is finite, so the Henselity of
$\bar \omega $ ensures that each finite extension $\Phi ^{\prime }$
of $\Phi _{\omega }$ contains as a subfield an inertial lift of 
$\widehat \Phi ^{\prime }$ over $\Phi _{\omega }$. At the same time, 
$\bar \omega $ is discrete, which shows that $\Phi ^{\prime }/\Phi 
_{\omega }$ is defectless if it is separable (see \cite{L}, Ch. XII, 
Sect. 6, Corollary~2). These facts make it easy to deduce from (2.5) 
and Galois theory that $\Gamma _{r}\Phi _{\omega }/\Phi _{\omega }$ 
is totally ramified and $[\Gamma _{r}K\colon K] = [\Gamma _{r}\Phi 
_{\omega }\colon \Phi _{\omega }] = [\Gamma _{r}\colon \Phi ] = p
^{\mu (r)}$. Therefore, $\Gamma _{r}/\Phi $ is totally ramified, i.e.
it possesses a primitive element $\theta $ whose minimal polynomial
$f _{\theta }(X)$ over $\Phi $ is Eisensteinian relative to $\omega $
(cf. \cite{FV}, Ch. 2, (3.6), and \cite{L}, Ch. XII, Sects. 2, 3 and 6).
Let $\theta _{0}$ be the free term of $f _{\theta }(X)$. As $\pi \in
\Phi $, $v(\pi ) \notin pv(K)$ and $\Gamma _{r}/\Phi $ is a Galois
extension, the conditions on $\theta $ guarantee that it is a
primitive element of $\Gamma _{r}K/K$ (and $\Gamma _{r}K _{v}/K
_{v}$), $p ^{\mu (r)}w(\theta ) = v(\theta _{0}) = \omega (\theta
_{0})$ and $v(\pi ) \in p ^{\mu (r)}w(\Gamma _{r}K)$, for any
valuation $w$ of $\Gamma _{r}K$ extending $v$. This implies $w$ is
unique, up-to an equivalence, and so completes the proof of Lemma
\ref{lemm2.3}.
\end{proof}
\par
\medskip
The conclusion of Lemma \ref{lemm2.3} need not be true in the 
mixed-characteristic setting. It has been established by Kurihara (cf. 
\cite{Ku}, Corollary~2) that there exists an HDV-field $(K, v)$ with 
char$(K) = 0$, $\widehat K$ imperfect and char$(\widehat K) = p > 0$, 
which does not admit a totally ramified cyclic extension of degree $p 
^{t}$, for any sufficiently large $t \in \mathbb N$ depending on $K$.

\medskip
\section{\bf Brauer $p$-dimensions in characteristic $p$}

\medskip
In this Section we consider index-exponent relations of
$p$-algebras over Henselian fields of characteristic $p$. First we
supplement Lemma \ref{lemm2.3} as follows:

\medskip
\begin{lemm}
\label{lemm3.1} Let $(K, v)$ be a valued field with {\rm char}$(K) =
p > 0$ and $v(K) \neq pv(K)$, and let $\tau (p)$ be defined as in
Theorem \ref{theo1.1}. Suppose that $L$ is a finite abelian extension
of $K$ in $K(p)$ satisfying the following conditions:
\par
{\rm (a)} $[L\colon K] = p ^{m}$ and $\mathcal{G}(L/K)$ has period
$p ^{m'}$ and rank $t$;
\par
{\rm (b)} $L$ has a unique, up-to an equivalence, valuation $v _{L}$
extending $v$, and the group $v _{L}(L)/v(K)$ is cyclic of order $p
^{m}$.
\par\noindent
Then there is $T \in d(K)$ with {\rm exp}$(T) = p ^{m'}$, possessing
a maximal subfield $K$-isomorphic to $L$, except, possibly, in case
$\tau (p) < \infty $ and $p ^{t-\tau (p)} \ge [\widehat K\colon
\widehat K ^{p}]$.
\end{lemm}

\medskip
\begin{proof}
It is clear from Galois theory and the structure of finite abelian
groups that $L = L _{1} \dots L _{t}$ and $[L\colon K] = \prod
_{j=1} ^{t} [L _{j}\colon K]$, for some cyclic extensions $L
_{j}/K$, $j = 1, \dots , t$. Take an element $\pi \in K$ so that
$v(\pi ) \in p ^{m}v _{L}(L)$, put $\pi _{0} = \pi $, and suppose
that there exist $\pi _{j} \in K ^{\ast }$, $j = 1, \dots , t$, and
$\mu \in \mathbb Z$ with $0 \le \mu \le t$, such that the cosets $v(\pi
_{i}) + pv(K)$, $i = 0, \dots , \mu $, are linearly independent over
$\mathbb F _{p}$, and in case $\mu < t$, $v(\pi _{u}) = 0$ and the
residue classes $\hat \pi _{u}$, $u = \mu + 1, \dots , t$, generate
an extension of $\widehat K ^{p}$ of degree $p ^{t-\mu }$ (this
assumption is admissible unless $\tau (p) \le t$ and $p ^{t-\tau (p)}
\ge [\widehat K\colon \widehat K ^{p}]$). Fix a generator $\lambda
_{j}$ of $\mathcal{G}(L _{j}/K)$, for $j = 1, \dots , t$, denote by
$T$ the $K$-algebra $\otimes _{j=1} ^{t} (L _{j}/K, \lambda _{j}, \pi
_{j})$, where $\otimes = \otimes _{K}$, and put $T ^{\prime } = T
\otimes _{K} K _{v}$. We show that $T \in d(K)$ (whence ind$(T) = p
^{m}$) and exp$(T) = p ^{m'}$. Clearly, $T ^{\prime } \cong \otimes
_{j=1} ^{t} (L _{j} ^{\prime }/K _{v}, \lambda _{j} ^{\prime }, \pi
_{j})$ over $K _{v}$, where $\otimes = \otimes _{K _{v}}$, $L _{j}
^{\prime } = L _{j}K _{v}$ and $\lambda _{j} ^{\prime }$ is the
unique $K _{v}$-automorphism of $L _{j} ^{\prime }$ extending
$\lambda _{j}$, for each $j$ (as in the proof of Lemma \ref{lemm2.3},
we identify $K _{v}$ with its $K$-isomorphic copy in $K _{\rm sep}$).
Therefore, it suffices for the proof of Lemma \ref{lemm3.1} to show
that $T ^{\prime } \in d(K _{v})$. Since, by the proof of Lemma
\ref{lemm2.3}, $K _{v}$ and $L ^{\prime } = LK _{v}$ are related as
in our lemma, this amounts to proving that $T \in d(K)$, for $(K, v)$
Henselian. Note that if $m = 1$, then our assertion is a special case
of \cite{Ch3}, Lemma~4.2. Henceforth, we assume that $m \ge 2$ and
view all value groups considered in the rest of the proof as
(ordered) subgroups of a fixed divisible hull of $v(K)$. Let $L _{0}$
be the degree $p$ extension of $K$ in $L _{t}$, and $R _{j} = L _{0}L
_{j}$, $j = 1, \dots , t$. Put $\rho _{t} = \lambda _{t} ^{p}$, and
when $t \ge 2$, denote by $\rho _{j}$ the unique $L
_{0}$-automorphism of $R _{j}$ extending $\lambda _{j}$, for $j = 1,
\dots , t - 1$. Then the centralizer $C$ of $L _{0}$ in $T$ is $L
_{0}$-isomorphic to $\otimes _{j=1} ^{t} (R _{j}/L _{0}, \rho _{j},
\pi _{j})$, where $\otimes = \otimes _{L _{0}}$; in particular,
deg$(C) = p ^{m-1}$. Using (2.1), Lemma \ref{lemm2.3} and this
result, one easily obtains that it is sufficient to prove that $T \in
d(K)$, under the extra hypothesis that $C \in d(L _{0})$.
\par
Let $w$ be the valuation of $C$ extending $v _{L _{0}}$, $\widehat
C$ its residue division ring, and for each $\xi \in C$ with $w(\xi )
= 0$, let $\widehat \xi \in \widehat C$ be the residue class of $\xi
$. It follows from the Ostrowski-Draxl theorem that $w(C)$ equals
the sum of $v(L)$ and the group generated by $[R _{i'}\colon L _{0}]
^{-1}v(\pi _{i'})$, $i ^{\prime } = 1, \dots , \mu $. Similarly, it
is proved that $\widehat C/\widehat K$ is a purely inseparable field
extension. Moreover, one sees that $\widehat C \neq \widehat K$ if
and only if $\mu < t$, and when this is the case, $[\widehat C\colon
\widehat K] = \prod _{u=\mu +1} ^{t} [R _{u}\colon L _{0}]$ and
\par\noindent
$\widehat C = \widehat K(\hat \eta _{\mu +1}, \dots , \hat \eta
_{t})$, where $\eta _{u} \in (R _{u}/L _{0}, \rho _{u}, \pi _{u})$ is
a root of $\pi _{u}$ of degree $[R _{u}\colon L _{0}]$ acting on $R
_{u}$ by conjugation as the automorphism $\rho _{u}$, for each index
$u$. In view of (2.1) and well-known general properties of purely
inseparable finite extensions (cf. \cite{L}, Ch. VII, Sect. 7), these
results show that $w(\eta _{t}) \notin pw(C)$, if $\mu = t$, and
$w(\eta _{t}) = 0$ and $\hat \eta _{t} \notin \widehat C ^{p}$,
otherwise. Observe now that there is a $K$-isomorphism $\bar \rho
_{t}$ of $C$ extending $\lambda _{t}$, such that $\bar \rho _{t}
^{p}(\bar c) = \eta _{t}\bar c\eta _{t} ^{-1}\colon \bar c \in C$,
and $\bar \rho _{t}(\eta _{t}) = \eta _{t}$. This implies $w(c) =
w(\bar \rho _{t}(c))$, for each $c \in C$, the products $c ^{\prime }
= \prod _{\kappa =0} ^{p-1} \bar \rho _{t} ^{\kappa }(c)$, $c \in C$,
have values $w(c ^{\prime }) \in pw(C)$, and $\hat c ^{\prime } \in
\widehat C ^{p}$, if $w(c) = 0$. Therefore, $c ^{\prime } \neq \eta
_{t}$, for any $c \in C$, so it follows from \cite{A2}, Ch. XI,
Theorems~11 and 12, that $T \in d(K)$. Let now $\Lambda $ be the
fixed field of the maximal subgroup of $\mathcal{G}(L/K)$ of period
$p$. Then \cite{P}, Sect. 15.1, Corollary~b, implies the class $p[D]
\in {\rm Br}(K)$ is represented by a crossed product of $\Lambda /K$,
defined similarly to $D$. As $\Lambda /K$ and $\pi $ are related like
$L/K$ and $\pi $, and $\mathcal{G}(\Lambda /K)$ is of period $p
^{m'-1}$, this enables one to prove inductively that exp$(D) = p
^{m'}$, as claimed.
\end{proof}

\medskip
\begin{coro}
\label{coro3.2}
Let $E$ be a field with {\rm char}$(E) = p > 0$ and $[E\colon E
^{p}] = p ^{\nu } < \infty $, and $F/E$ a finitely-generated
extension of transcendency degree $n > 0$. Then $n + \nu - 1 \le
{\rm Brd}_{p}(F) \le {\rm abrd}_{p}(F) \le n + \nu $, and when $n +
\nu \ge 2$,
\par\noindent
$(p ^{t}, p ^{s})\colon t, s \in \mathbb N, s \le t \le (n + \nu -
1)s$, are index-exponent pairs over $F$.
\end{coro}

\medskip
\begin{proof}
We have $n + \nu - 1 \le {\rm Brd}_{p}(F) \le {\rm abrd}_{p}(F) \le
n + \nu $, by \cite{Ch3}, Theorem~2.1 (c). Note also that $F$ has a
valuation $v$ trivial on $E$, such that $v(F) = \mathbb Z ^{n}$ and
$\widehat F$ is a finite extension of $E$ (see, e.g. \cite{Ch3},
(4.1)). Therefore, $[\widehat F\colon \widehat F ^{p}] = p ^{\nu }$
(cf. \cite{L}, Ch. VII, Sect. 7) and $v(F)/pv(F)$ is of order $p
^{n}$, which makes it easy to deduce the concluding assertion of
Corollary \ref{coro3.2} from Lemma \ref{lemm3.1}.
\end{proof}

\medskip
\begin{rema}
\label{rema3.3} It is known \cite{PY}, (3.19) (see also \cite{JW},
Corollary~6.10) that if $(K, v)$ is a Henselian field and $T \in
d(K)$ is a tame $K$-algebra, in the sense of \cite{PY} or \cite{JW},
then the period per$(T/K)$ of the group $v(T)/v(K)$ divides
exp$(T)$. At the same time, by Lemma \ref{lemm3.1} with its proof,
if char$(K) = p > 0$ and $v(K)/pv(K)$ is infinite, then there are
$T _{n} \in d(K)$, $n \in \mathbb N$, such that ind$(T _{n}) = {\it
per}(T _{n}/K) = p ^{n}$, exp$(T _{n}) = p$ and $T _{n}/K$ is
defectless, for each $n$.
\end{rema}

\medskip
Next we describe index-exponent $p$-primary pairs over some 
maximally complete fields of characteristic $p$, including those with 
perfect residue fields.

\medskip
\begin{prop}
\label{prop3.4} Let $(K, v)$ be a valued field of characteristic $p
> 0$. Suppose that $v(K)/pv(K)$ is infinite or $[\widehat K\colon
\widehat K ^{p}] = \infty $, where $\widehat K ^{p} = \{\hat a
^{p}\colon \ \hat a \in \widehat K\}$. Then $(p ^{k}, p ^{n})\colon
k, n \in \mathbb N, n \le k$, are index-exponent $K$-pairs.
\end{prop}

\medskip
\begin{proof}
Lemma \ref{lemm3.1}, \cite{Ch5}, Remark~4.3, and our assumptions show
that there are tensor products $D _{n} \in d(K)$, $n \in \mathbb N$,
of degree $p$ cyclic $K$-algebras with exp$(D _{n}) = p$ and ind$(D
_{n}) = p ^{n}$, for each $n$. Hence, by \cite{Ch4}, Lemma~5.2, it
suffices to prove that $(p ^{n}, p ^{n})$, $n \in \mathbb N$, are
index-exponent $K$-pairs. By Witt's lemma (cf. \cite{Dr1}, Sect. 15,
Lemma~2), each cyclic extension $L$ of $K$ in $K(p)$ lies in $I(L
^{\prime }/K)$, for some $\mathbb Z _{p}$-extension $L ^{\prime }$ of
$K$ in $K(p)$. Fix a topological generator $\sigma $ of
$\mathcal{G}(L ^{\prime }/K)$, and for any $n \in \mathbb N$, let $L
_{n}$ be the extension of $K$ in $L ^{\prime }$ of degree $p ^{n}$,
and $\sigma _{n}$ the automorphism of $L _{n}$ induced by $\sigma $.
Clearly, $L _{n}/K$ is cyclic and $\sigma _{n}$ generates
$\mathcal{G}(L _{n}/K)$. Choosing $L ^{\prime }$ so that $(L _{1}/K,
\sigma _{1}, c) \cong D _{1}$, for some $c \in K ^{\ast }$, one gets
ind$(\Delta _{n}) = {\rm exp}(\Delta _{n}) = p ^{n}$ from \cite{P},
Sect. 15.1, Corollary~b, for the cyclic $K$-algebras $\Delta _{n} =
(L _{n}/K, \sigma _{n}, c)$, $n \in \mathbb N$, which completes our
proof.
\end{proof}

\medskip
\begin{prop}
\label{prop3.5} Let $(K, v)$ be a maximally complete field with {\rm
char}$(K) = p > 0$, $v(K) \neq pv(K)$ and $[K\colon K ^{p}] = p
^{n}$, for some $n \in \mathbb N$, and let $G _{p}$ be a Sylow
pro-$p$-subgroup of $\mathcal{G}(\widehat K _{\rm sep}/\widehat K)$.
Then $n - 1 \le {\rm Brd}_{p}(K) \le n$. Moreover, the following holds
when $\widehat K$ is perfect:
\par
{\rm (a)} {\rm Brd}$_{p}(K) = n - 1$ if and only if $n > r _{p}(\widehat
K)$;
\par
{\rm (b)} $(p ^{k}, p ^{s})\colon k, s \in \mathbb N$, $s \le k \le
{\rm Brd}_{p}(K)s$, are index-exponent $K$-pairs.
\par
{\rm (c)} {\rm abrd}$_{p}(K) = n - 1$ if and only if either $G _{p} =
\{1\}$ or $n \ge 2$ and $G _{p} \cong \mathbb Z _{p}$.
\end{prop}

\medskip
\begin{proof}
Our assumptions show that $[K\colon K ^{p}] = [\widehat K\colon
\widehat K ^{p}]e(K/K ^{p})$ (cf. \cite{W}, Theorem~31.21), so it
follows from Lemma \ref{lemm3.1} and Albert's theory of $p$-algebras
\cite{A2}, Ch. VII, Theorem~28, that $n - 1 \le {\rm Brd}_{p}(K) \le
n$, as claimed. In the rest of the proof, we suppose that $\widehat
K$ is perfect. First we consider the case of $r _{p}(\widehat K)
\ge n$. Then one gets from Galois theory and Witt's lemma that
$\mathbb Z _{p} ^{n}$ is realizable as a Galois group over $\widehat
K$. Hence, by \cite{TW}, Theorem A.24, there is a Galois extension $U
_{n}$ of $K$ in $K _{\rm ur}$ with $\mathcal{G}(U _{n}/K) \cong
\mathbb Z _{p} ^{n}$. This implies each finite abelian $p$-group $H$
of rank $\le n$ is isomorphic to $\mathcal{G}(U _{H}/K)$, for some
Galois extension $U _{H}$ of $K$ in $U _{n}$. Observing also
that $v(K)/pv(K)$ has order $p ^{n}$, and using \cite{JW},
Example~4.3, one proves the existence of an NSR-algebra $N _{H} \in
d(K)$ with a maximal subfield $U _{H} ^{\prime } \cong U _{H}$ over 
$K$. Therefore, exp$(N _{H}) = {\rm per}(H)$ and ind$(N _{H}) = [U 
_{H}\colon K]$, so Brd$_{p}(K) = n$, which reduces the rest of our 
proof to the case of $n > r _{p}(\widehat K)$. Note that $(L ^{\prime 
}, v _{L'})$ is maximally complete and $[L ^{\prime }\colon L 
^{\prime p}] = p ^{n}$ whenever $L ^{\prime }/K$ is a finite 
extension (cf. \cite{W}, Theorem~31.22, and \cite{L}, Ch. VII, Sect. 
7). This enables one to deduce from \cite{AJ}, Theorem~3.3, by the 
method of proving \cite{Ch5}, (5.5), that for each $D _{e} \in d(K)$ 
with exp$(D _{e}) = p ^{e}$, where $e \in \mathbb N$, $[D _{e}] \in 
{\rm Br}(K _{e}/K)$, for some purely inseparable extension $K _{e}/K$ 
such that $[K _{e}\colon K] \mid p ^{(n-1)e}$. In view of (1.1) (b), 
the obtained result yields ind$(D _{e}) \mid p ^{(n-1)e}$ and 
Brd$_{p}(K) = n - 1$, so Proposition \ref{prop3.5} (a) is proved. 
Applying Lemmas \ref{lemm2.3} and \ref{lemm3.1}, one concludes that 
$(p ^{t}, p ^{m})$, $t, m \in \mathbb N$, $0 < m \le t \le (n - 1)m$, 
are index-exponent $K$-pairs, which reduces Proposition \ref{prop3.5} 
(b) to a consequence of Proposition \ref{prop3.5} (a). It remains for 
us to prove Proposition \ref{prop3.5} (c). Clearly, if $G _{p} = 
\{1\}$, then $r _{p}(\widehat L) = 0$, for every $L \in {\rm Fe}(K)$. 
At the same time, it follows from Galois cohomology and 
Nielsen-Schreier's formula for open subgroups of free pro-$p$-groups 
(cf. \cite{S2}, Ch. I, 3.3, 4.2; Ch. II, 2.2) that if $G _{p}$ is not 
procyclic, then $r _{p}(K _{1}) \ge n$, for some finite extension $K 
_{1}$ of $K$ in $K _{\rm ur}$. Note finally that if $G _{p}$ has 
rank $1$ as a pro-$p$-group, then its open subgroups are isomorphic 
to $\mathbb Z _{p}$, which implies $r _{p}(L) \le 1$, $L \in {\rm 
Fe}(K)$. As $(L, v _{L})$ is maximally complete and $[L\colon K] = p 
^{n}$, these facts give us the possibility to deduce Proposition 
\ref{prop3.5} (c) from Proposition \ref{prop3.5} (a).
\end{proof}

\medskip
We are now prepared to generalize Proposition \ref{prop1.2} as follows.
\medskip
\begin{coro}
\label{coro3.6} Let $(K, v)$ be a maximally complete field with {\rm
char}$(K) = p > 0$ and $\tau (p)$ defined as in Theorem
\ref{theo2.1}. Suppose further that $\widehat K$ is complete with
respect to a discrete valuation $\omega $ with a quasifinite residue
field $\tilde k$. Then:
\par
{\rm (a)} {\rm Brd}$_{p}(K) = \infty $ if and only if $\tau (p) =
\infty $; when this holds, $(p ^{k}, p ^{n})$ is an index-exponent
pair over $K$, for any $k, n \in \mathbb N$ with $k \ge n$;
\par
{\rm (b)} {\rm Brd}$_{p}(K) = \tau (p)$, provided that $\tau (p) <
\infty $; in this case, $(p ^{k}, p ^{n})$ is an index-exponent
$K$-pair, where $k, n \in \mathbb N$, if and only if $n \le k \le
n\tau (p)$.
\end{coro}

\medskip
\begin{proof}
It is known (cf. \cite{Efr2}, Sect. 5.2) that $K$ has a valuation
$\varphi $ (a refinement of $v$), such that $\varphi (K) = v(K)
\oplus \omega (\widehat K)$, $\omega (\widehat K)$ is an isolated
subgroup of $\varphi (K)$, $v$ and $\omega $ are canonically induced
by $\varphi $ and $\omega (\widehat K)$ on $K$ and $\widehat K$,
respectively, and $\widehat K _{\varphi } \cong \tilde k$, where
$\widehat K _{\varphi }$ is the residue field of $(K, \varphi )$.
Observing that, by theorems of Krull and Hasse-Schmidt-MacLane
(cf. \cite{Efr2}, Theorems~12.2.3, 18.4.1, and \cite{W},
Theorem~31.24 and page 483), $(\widehat K, \omega )$ is maximally
complete and $(K, \varphi )$ possesses an immediate extension $(K 
^{\prime }, \varphi ^{\prime })$ which is a maximally complete field, 
one obtains that $(K ^{\prime }, \varphi ^{\prime }) = (K, \varphi 
)$. As $r _{p}(\tilde k) = 1$ and $\tilde k$ is perfect, Corollary 
\ref{coro3.6} can now be deduced from Propositions \ref{prop3.4} and 
\ref{prop3.5}.
\end{proof}

\medskip
When $(K, v)$ is a Henselian field, such that char$(K) = p > 0$,
$v(K)$ is a non-Archimedean group, $v(K)/pv(K)$ is finite and
$[\widehat K\colon \widehat K ^{p}] = p ^{\nu } < \infty $, there
is, generally, no formula for Brd$_{p}(K)$ involving only invariants
of $\widehat K$ and $v(K)$. This is illustrated below in the case of
$v(K) = \mathbb Z ^{t}$, for any integer $t \ge 2$.

\medskip
\begin{exam}
Let $Y _{0}$ be a field with char$(Y _{0}) = p$ and $[Y _{0}\colon Y
_{0} ^{p}] = p ^{\nu } < \infty $, and let $Y _{t} = Y _{0}((T _{1})) 
\dots ((T _{t}))$ be the iterated formal Laurent power series field in $t$
variables over $Y _{0}$. Denote by $w _{t}$ the natural $\mathbb Z
^{t}$-valued valuation of $Y _{t}$ trivial on $Y _{0}$. It is known
(see \cite{BK}, page 181 and further references there) that there 
exist elements $X _{n} \in Y _{t-1}$, $n \in \mathbb N$, algebraically
independent over the field $Y _{t-2}(T _{t-1})$, where $Y _{t-2}
= Y _{0}((T _{1})) \dots ((T _{t-2}))$ in the case of $t \ge 3$. Put
$F _{n} = Y _{t-2}(T _{t-1}, X _{1}, \dots X _{n})$, for each $n \in
\mathbb N$, $F _{\infty } = \cup _{n=1} ^{\infty } F _{n}$, and
$\mathbb N _{\infty } = \mathbb N \cup \{\infty \}$. For any $n \in
\mathbb N _{\infty }$, denote by $F _{n} ^{\prime }$ the separable
closure of $F _{n}$ in $Y _{t-1}$, and by $v _{n}$ the valuation of
the field $K _{n} = F _{n} ^{\prime }((T _{t}))$ induced by $w
_{t}$. It is easily verified that $(K _{n}, v _{n})$ is Henselian, $v 
_{n}(K _{n}) = \mathbb Z ^{t}$ and $\widehat K _{n} = Y 
_{0}$, for each index $n$. Note also that $[F _{\infty } ^{\prime
}\colon F _{\infty } ^{\prime p}] = \infty $, so Proposition
\ref{prop3.4}, applied to the valuation of $K _{n}$ induced by the
natural discrete valuation of $Y _{t}$ trivial on $Y _{t-1}$, yields
Brd$_{p}(K _{\infty }) = \infty $. When $n \in \mathbb N$, we have
$[K _{n}\colon K _{n} ^{p}] = p ^{\nu + t+n} = p[F _{n} ^{\prime
}\colon F _{n} ^{\prime p}]$, which enables one to deduce from Lemma
\ref{lemm3.1}, \cite{Ch3}, Lemma~4.1, and \cite{A2}, Ch. VII,
Theorem~28 (see also \cite{L}, Ch. VII, Sect. 7) that $\nu + t + n -
1 \le {\rm Brd}_{p}(K _{n}) \le \nu + n + t$.
\end{exam}

\medskip
\section{\bf Brauer $p$-dimensions of Henselian fields with
$p$-quasilocal residue fields}

\medskip
Let $(K, v)$ be a Henselian field with $\widehat K$ $p$-quasilocal
and $r _{p}(\widehat K) > 0$. Then Brd$_{p}(\widehat K) \le 1$, so
Theorem \ref{theo2.1} yields Brd$_{p}(K) = \infty $ if and only if
$m _{p} = \infty $ or $\tau (p) = \infty $ and $\varepsilon _{p} \in
\widehat K$. When Brd$_{p}(K) = \infty $, index-exponent $p$-primary
$K$-pairs are described by Corollary \ref{coro2.2} (and the 
Pythagorean property of formally real $2$-quasilocal fields, see 
\cite{Ch2}, I, Lemma~3.5). The main result of this Section concerns
the case of Brd$_{p}(K) < \infty $ and can be stated as follows:

\medskip
\begin{theo}
\label{theo4.1} Let $(K, v)$ be a Henselian field with {\rm
Brd}$_{p}(K) < \infty $, for some $p \in \mathbb P$, and set
$\varepsilon _{p}$, $\tau (p)$ and $m _{p}$ as in Theorem
\ref{theo2.1}. Suppose that $\widehat K$ is $p$-quasilocal, $p \neq
{\rm char}(\widehat K)$ and $m _{p} > 0$. Then:
\par
{\rm (a)} ${\rm Brd}_{p}(K) = u _{p}$, where $u _{p} =
[(\tau (p) + m _{p})/2]$, if $\varepsilon _{p} \in \widehat K$ and
$\widehat K$ is a nonreal field; $u _{p} = m _{p}$, if $\varepsilon _{p}
\notin \widehat K$;
\par
{\rm (b)} {\rm Br}$(K) _{2}$ is a group of period $2$ and {\rm
Brd}$_{2}(K) = 1 + [\tau (2)/2]$, provided that $\widehat K$ is
formally real and $p = 2$.
\end{theo}
\par
\medskip
Before proving Theorem \ref{theo4.1}, note that it yields Brd$_{p}(K)
= \tau (p)$ whenever $r _{p}(\widehat K) = \infty $. This holds in 
all presently known cases where $\widehat K$ is $p$-quasilocal and 
Br$(\widehat K) _{p}$ does not embed in $\mathbb Q/\mathbb Z$ or, 
equivalently, in the quasicyclic $p$-group $\mathbb Z(p ^{\infty })$ 
(see \cite{VdBS}, the end of Sect.~3, \cite{Ch6}, Theorem~1.2, and 
e.g., \cite{M}, \cite{T}).

\medskip
{\it Proof of Theorem \ref{theo4.1}}. Suppose first that $\widehat K$
is formally real and $p = 2$. Then, by \cite{Ch2}, I, Lemma~3.5,
$\widehat K$ is Pythagorean, $\widehat K(2) = \widehat K(\sqrt{-1})$
and Br$(\widehat K) _{2}$ is of order $2$. Therefore, $r
_{2}(\widehat K) = 1$ and $r _{2}(\widehat K(\sqrt{-1})) = 0$, so it
follows from the Merkur'ev-Suslin theorem \cite{MS}, (16.1), that
Br$(\widehat K(\sqrt{-1})) _{2} = \{0\}$. Note further that $K$ is
Pythagorean, which implies $2{\rm Br}(K) = \{0\}$ (cf. \cite{La},
Theorem~3.16, and \cite{Efr1}, Theorem~3.1). These observations and
\cite{Ch5}, Corollary~5.5, prove Theorem \ref{theo4.1} (b). We turn
to the proof of Theorem \ref{theo4.1} (a), so we assume that $p > 2$
or $\widehat K$ is a nonreal field. Then Br$(\widehat K) _{p}$ is a
divisible group, by \cite{Ch2}, I, Theorem~3.1. Our argument also
relies on the following results concerning inertial algebras $I \in 
d(K)$ with $[I] \in {\rm Br}(K) _{p}$, and inertial extensions $U$ of 
$K$ in $K(p)$:
\par
\medskip
(4.1) (a) ind$(I) = {\rm exp}(I)$ and $I$ is a cyclic $K$-algebra;
\par
(b) $[I] \in {\rm Br}(U/K)$ if and only if ind$(I) \mid [U\colon
K]$; $U$ is embeddable in $I$ as a $K$-subalgebra if and only if
$[U\colon K] \mid {\rm ind}(I)$;
\par
(c) ind$(I \otimes _{K} I ^{\prime })$ equals ind$(I)$ or ind$(I
^{\prime })$, if $I ^{\prime } \in d(K)$, $I ^{\prime }/K$ is NSR,
and $[I ^{\prime }] \in {\rm Br}(K) _{p}$.
\par
\medskip\noindent
Statements (4.1) can be deduced from (1.3), (2.3) (b) and \cite{JW},
Theorems~3.1 and 5.15. They imply in conjunction with \cite{Ch5},
Lemma~4.1, that ind$(W) \mid {\rm exp}(W) ^{m _{p}}$, for each $W
\in d(K)$ inertially split over $K$. At the same time, it follows
from \cite{Ch3}, (3.3) and (3.6), and \cite{Mo}, Theorem~1 (see also
\cite{JW}, Example~4.3), that there is an NSR-algebra $W ^{\prime }
\in d(K)$ with ind$(W ^{\prime }) = p ^{m _{p}}$ and exp$(W ^{\prime
}) = p$. Observe now that, by (2.3) (c), Br$(K) _{p} \subseteq {\rm
Br}(K _{\rm ur}/K)$ in case $\varepsilon _{p} \notin \widehat K$ or
$\tau (p) = 1$. In view of (4.1) and \cite{JW}, Theorem~4.4 and
Lemma~5.14, this yields Brd$_{p}(K) = m _{p}$, so it remains for us
to prove Theorem \ref{theo4.1}, under the extra hypothesis that
$\varepsilon _{p} \in \widehat K$ and $\tau (p) \ge 2$. It is easily
obtained from \cite{Mo}, Theorem~1, and \cite{Ch5}, Lemmas~4.1 and
4.2, that there exists $\Delta \in d(K)$ with exp$(\Delta ) = p$ and
ind$(\Delta ) = p ^{\mu (p)}$, where $\mu (p) = [(m _{p} + \tau
(p))/2]$. This means that Brd$_{p}(K) \ge \mu (p)$, so we have to
prove that Brd$_{p}(K) \le \mu (p)$. Note first that $2 \le m _{p}$, 
provided Br$(\widehat K) _{p} \neq \{0\}$. Assuming the opposite and 
taking into account that $\varepsilon _{p} \in \widehat K$, one 
obtains from the other conditions on $\widehat K$ that it is a 
nonreal field with $r _{p}(\widehat K) = 1$. Hence, by \cite{Wh}, 
Theorem~2, $\widehat K(p)/\widehat K$ is a $\mathbb Z 
_{p}$-extension. In view of \cite{MS}, (11.5) and (16.1), and Galois
cohomology (cf. \cite{S2}, Ch. I, 4.2), this requires that 
Br$(\widehat K) _{p} = \{0\}$. As $\tau (p) \ge 2$, the obtained
contradiction proves that $r _{p}(\widehat K) \ge m _{p} \ge 2$, as 
claimed. Now take an algebra $D \in d(K)$ so that exp$(D) = p ^{n}$, 
for some $n \in \mathbb N$, attach $S$, $V$ and $T \in d(K)$ to $D$ 
as in (2.3) (a), and fix $\Theta \in d(K)$ so that $[\Theta ] = [V 
\otimes _{K} T]$. To prove that ind$(D) \mid p ^{n\mu (p)}$ we need 
the following statements:
\par
\medskip
(4.2) (a) If $n = 1$, then $S$, $V$ and $T$ can be chosen so that $V
\otimes _{K} T = \Theta $, and $S = K$ or $V = K$.
\par
(b) If $n \ge 2$, then there is a totally ramified extension $Y$ of 
$K$ in $K(p)$, such that $[Y\colon K] \mid p ^{\mu (p)}$ and either 
exp$(D _{Y}) \mid p ^{n-1}$, or exp$(D _{Y}) = {\rm exp}(S _{Y}) = p 
^{n}$, $[Y\colon K] \mid p ^{[\tau (p)/2]}$ and exp$(V _{Y} \otimes 
_{Y} T _{Y}) \mid p ^{n-1}$, where $S _{Y}, V _{Y}, T _{Y} \in d(Y)$ 
are attached in accordance with (2.3) (a) to a representative $D _{Y} 
\in d(Y)$ of $[D \otimes _{K} Y]$.
\par
\medskip\noindent
Statement (4.2) (a) can be deduced from (4.1), \cite{Ch5}, (4.7),
and well-known properties of cyclic algebras (cf. \cite{P}, Sect.
15.1, Proposition~b). Since $m _{p} \ge 2$, (4.2) (a) implies the 
assertion of Theorem \ref{theo4.1} (a) in the case of $n = 1$, so we 
assume further that $n \ge 2$. The conclusion of (4.2) (b) is obvious, if
exp$(\Theta ) \mid p ^{n-1}$ (one may put $Y = K$). Therefore, by
(2.4) (c), it suffices to prove (4.2) (b) under the hypothesis that
exp$(\Theta ) = p ^{n}$. Take $D _{n-1} \in d(K)$ so that $[D
_{n-1}] = p ^{n-1}[D]$ and attach to it a triple $S _{n-1}$, $V
_{n-1}$, $T _{n-1} \in d(K)$ in agreement with (4.2) (a). Then $V
_{n-1} \otimes _{K} T _{n-1}$ contains as a maximal subfield an
abelian and totally ramified extension $Y$ of $K$. Observing that $[V 
_{n-1} \otimes _{K} T _{n-1}] \in {\rm Br}(Y/K)$, identifying $Y$
with its $K$-isomorphic copy in $K(p)$, and using (2.4) (a) and (1.1) 
(a), one sees that it has the properties required by (4.2) (b).
\par
We continue with the proof of Theorem \ref{theo4.1} (a). In view of 
(2.2) and (4.2) (a), a standard inductive argument allows us to 
proceed under the extra hypothesis that ind$(D ^{\prime }) \mid {\rm 
exp}(D ^{\prime }) ^{\mu (p)}$, for each $D ^{\prime } \in d(K 
^{\prime })$ with exp$(D ^{\prime }) \mid p ^{n-1}$, where $K 
^{\prime }/K$ is an arbitrary totally ramified finite extension. It is 
known (cf. \cite{JW}, Corollary~6.8) that if $J, J ^{\prime } \in 
d(K)$, $J/K$ is inertial and $[J ^{\prime }] = [J \otimes _{K} \Theta 
]$, then $v(J ^{\prime }) = v(\Theta )$, $Z(\widehat J ^{\prime }) = 
Z(\widehat \Theta )$ and $[\widehat J ^{\prime }] = [\widehat J 
\otimes _{\widehat K} \widehat \Theta ] \in {\rm Br}(Z(\widehat 
\Theta ))$. Note also that the period of the group $v(J ^{\prime 
})/v(K)$ divides exp$(J ^{\prime })$  (see Remark \ref{rema3.3}). At 
the same time, by \cite{Ch2}, I, Theorem~4.1, the scalar extension 
map Br$(\widehat K) \to {\rm Br}(Z(\widehat \Theta ))$ induces a 
surjective homomorphism Br$(\widehat K) _{p} \to {\rm Br}(Z(\widehat 
\Theta )) _{p}$. As Brd$_{p}(\widehat K) \le 1$ and $m _{p} \ge 2$, 
these results, combined with (1.3), (4.1) (a), (b), the 
Ostrowski-Draxl theorem, and the inductive hypothesis, prove the 
following:
\par
\medskip
(4.3) (a) If exp$(\Theta ) \mid p ^{n-1}$, then ind$(D) \mid p.{\rm
ind}(S _{0} \otimes _{K} V \otimes _{K} T)$, for some $S _{0} \in
d(K)$ inertial over $K$ with exp$(S _{0}) \mid p ^{n-1}$;
\par
(b) If exp$(\Theta ) \mid p ^{n-1}$ and ind$(D) > {\rm ind}(I
\otimes _{K} V \otimes _{K} T)$ whenever $I \in d(K)$, $[I] \in {\rm
IBr}(K)$ and exp$(I) \mid p ^{n-1}$, then $[Z(\widehat D)\colon
\widehat K] = p ^{k}$ and $[\widehat D\colon Z(\widehat D)] = p
^{2n-2k}$, for some $k \in \mathbb Z$ with $0 \le k < n$; hence,
ind$(D) ^{2} \mid p ^{2n}e(\Theta /K) \mid p ^{2n}{\rm exp}(\Theta )
^{\tau (p)}$, which yields ind$(D) ^{2} \mid p ^{2n+(n-1)\tau (p)}
\mid p ^{m _{p}n+(n-1)\tau (p)}$.
\par
\medskip\noindent
Now fix an extension $Y/K$ and $Y$-algebras $D _{Y}$, $S _{Y}$, $V
_{Y}$, $T _{Y}$ as in (4.2) (b), and take $\Theta _{Y} \in d(Y)$ so
that $[\Theta _{Y}] = [V _{Y} \otimes _{Y} T _{Y}]$. Observing that,
by (1.1) (b), ind$(D) \mid {\rm ind}(D _{Y})[Y\colon K]$, and
applying (4.3) in case exp$(D _{Y}) = p ^{n}$ to $D _{Y}$, $V _{Y}$,
$T _{Y}$ and $\Theta _{Y}$, instead of $D$, $V$, $T$ and $\Theta $,
respectively, one concludes that ind$(D) ^{2} \mid p ^{n(m _{p}+\tau
(p))}$. Theorem \ref{theo4.1} is proved.
\par
\medskip
Theorem \ref{theo4.1} (a) retains its validity, if $(K, v)$ is a Henselian
field, such that $\tau (p) < \infty $, $r _{p}(\widehat K) = 0$ and 
$\mu _{p}(\widehat K) \neq \{1\}$. Then it follows from \cite{MS}, 
(16.1), that Brd$_{p}(\widehat K) = 0$, so Theorem \ref{theo2.1} (a) 
implies Brd$_{p}(K) = [\tau (p)/2]$.
\par
\medskip
\begin{rema}
\label{rema4.2}
Let $(K, v)$ be a Henselian field with $\widehat K$ formally real and
$2$-quasilocal. Then the symbol $K$-algebra $D ^{\prime } = A
_{-1}(-1, -1; K)$ lies in $d(K)$, and it follows from \cite{Ch5},
Lemma~4.2, that if $\tau (2) \ge 2$, then there exist $D _{n} \in
d(K)$, $n = 1, \dots , [\tau (2)/2]$, totally ramified over $K$ with
exp$(D _{n}) = 2$ and ind$(D _{n}) = 2 ^{n}$, for each $n$.
As $D ^{\prime }/K$ is inertial, this implies together with
\cite{Mo}, Theorem~1, that $D ^{\prime } \otimes _{K} D _{n} \in
d(K)$ (and ind$(D ^{\prime } \otimes _{K} D _{n}) = 2 ^{n+1}$),
$n = 1, \dots , [\tau (2)/2]$. In view of (2.3) (b) and Theorem
\ref{theo4.1} (b), these facts prove that if $0 \le \tau (2) <
\infty $, then $(1, 1)$ and $(2 ^{n}, 2)$, $n = 1, \dots , 1 + [\tau
(2)/2]$, are all index-exponent $2$-primary $K$-pairs.
\end{rema}
\par
\medskip
\begin{coro}
\label{coro4.3}
Let $K _{m}$ be an $m$-dimensional local field with a quasifinite
$m$-th residue field $K _{0}$, for some $m \in \mathbb N$. Suppose
that $p \in \mathbb P$ is different from {\rm char}$(K _{0})$, and
$\varepsilon _{p}$ is a primitive $p$-th root of unity in $K _{0,{\rm
sep}}$. Then {\rm Brd}$_{p}(K _{m}) = [(1 + m)/2]$, if $\varepsilon
_{p} \in K _{0}$; {\rm Brd}$_{p}(K _{m}) = 1$, otherwise.
\end{coro}
\par
\medskip
\begin{proof}
This is in fact a special case of Theorem \ref{theo4.1}, since our
assumptions imply the existence of a Henselian $\mathbb Z
^{m}$-valued valuation on $K _{m}$ with $\widehat K _{m} = K _{0}$.
\end{proof}
\par
\medskip
When $\varepsilon _{p} \in K _{0}$, the equality Brd$_{p}(K _{m}) =
[(1 + m)/2]$ can also be obtained from \cite{Ch3}, Lemma~4.1, and Khalin's
formula for the number of isomorphism classes of $K _{m}$-algebras $D
_{p,k} \in d(K _{m})$ with exp$(D _{p,k}) = p$ and ind$(D _{p,k}) = p
^{k}$, for a fixed $k \in \mathbb N$ (Khalin's formula has been
deduced in \cite{Kh1}, under the hypothesis that $K _{0}$ is finite,
but it clearly holds in the setting of Corollary \ref{coro4.3} as well).
\par
\medskip
\begin{prop}
\label{prop4.4}
Let $K _{m}$ be an $m$-dimensional local field with {\rm char}$(K
_{m}) = 0$, $K _{0}$ finite and {\rm char}$(K _{0}) = p$. Then $m
- 1 \le {\rm abrd}_{p}(K _{m}) \le m$. Moreover, {\rm Brd}$_{p}(K
_{m}) \ge m - 1$ unless $m \ge 4$, {\rm char}$(K _{1}) = 0$ and $r
_{p}(K _{1}) < m - 1$, where $K _{1}$ is the last but one residue
field of $K _{m}$.
\end{prop}
\par
\medskip
\begin{proof}
Note that if $m = 1$, then Brd$_{p}(K _{m}) = {\rm abrd}_{p}(K _{m}) 
= 1$ (cf. \cite{S1}, Ch. XIII, Sect. 3), which proves our assertions. We 
assume further that $m \ge 2$. It is well-known that finite 
extensions of $K _{m}$ are $m$-dimensional local fields, so the 
equality abrd$_{p}(K _{m}) \le m$ reduces to a consequence of 
\cite{Ch4}, Lemma~4.1, and the Corollary to \cite{Kh2}, Theorem~2. 
To prove the other inequalities stated in Proposition \ref{prop4.4}, 
we consider the $i$-th residue field $K _{m-i}$ of $K _{m}$, where 
$i \ge 0$ is the maximal integer for which char$(K _{m-i}) = 0$. 
Clearly, if $i > 0$, then $K _{m}$ has a $\mathbb Z ^{i}$-valued 
Henselian valuation $v _{i}$ with a residue field $K _{m-i}$. When 
$i = m - 1$, Theorem \ref{theo4.1}, applied to $(K _{m}, v _{i})$, 
gives a formula for Brd$_{p}(K _{m})$, which indicates that 
Brd$_{p}(K _{m}) \le m - 1$ and equality holds if and only if $r 
_{p}(K _{1}) \ge m - 1$. This, combined with \cite{S2}, Ch. II, 
Theorems~3 and 4 (applied to finite extensions of $K _{1}$), proves 
that abrd$_{p}(K _{m}) = m - 1$. Thus it follows that Brd$_{p}(K 
_{m}) = m - 1$ in case $m \le 3$. It remains to be seen that 
Brd$_{p}(K _{m}) \ge m - 1$, provided that $i < m - 1$. Then $K 
_{m-i'}$, $i ^{\prime } = i, i + 1$, is an $(m - i ^{\prime 
})$-dimensional local field with last residue field $K _{0}$; in 
particular, $K _{m-i'}$ is complete with respect to a discrete 
valuation $\omega _{m-i'}$ whose residue field is $K _{m-i'-1}$. In 
view of Lemma \ref{lemm2.3} and Proposition \ref{prop3.5}, this 
means that $r _{p}(K _{m-i-1}) = \infty $, and in the case where $i 
< m - 2$, Brd$_{p}(K _{m-i-1}) = m - i - 2$. More precisely, there 
exist $D _{0} \in d(K _{m-i-1})$, defined as in the proof of Lemma 
\ref{lemm3.1} when $i < m - 2$ (and equal to $K$, if $i = m - 2$), 
and totally ramified Galois extensions $M _{n} ^{\prime }/K 
_{m-i-1}$, $n \in \mathbb N$, relative to $\omega _{m-i-1}$, such 
that deg$(D _{0}) = e(D _{0}/K _{m-i-1}) = p ^{m-i-2}$, $[D _{0}] 
\in $ $_{p}{\rm Br}(K _{m-i-1})$, $\widehat D _{0}$ is a field with 
$\widehat D _{0} ^{p} \subseteq \widehat K$, and for each index $n$, 
$D _{0} \otimes _{K_{m-i-1}} M _{n} ^{\prime } \in d(M _{n} ^{\prime 
})$ and $\mathcal{G}(M _{n} ^{\prime }/K _{m-i-1})$ is elementary 
abelian of order $p ^{n}$. Let $D$ and $M _{n}$ be inertial lifts 
over $K _{m-i}$ (relative to $\omega _{m-i}$) of $D _{0}$ and $M 
_{n} ^{\prime }$, respectively. Then $M _{n}/K _{m-i}$ are inertial 
Galois extensions, $\mathcal{G}(M _{n}/K _{m-i}) \cong \mathcal{G}(M 
_{n} ^{\prime }/K _{m-i-1})$ and $D \otimes _{K _{m-i}} M _{n}$ lies 
in $d(M _{n})$, for every $n \in \mathbb N$. This enables one to 
deduce (in the spirit of the proof of \cite{Ch5}, Proposition~6.3) 
from \cite{JW}, Example~4.3 (or \cite{Ch4}, (3.6) (a)), and 
\cite{Mo}, Theorem~1, that there exists $T \in d(K _{m-i})$ with 
deg$(T) = p$, $T/K _{m-i}$ NSR relative to $\omega _{m-i}$, and 
$\Sigma \in d(K _{m-i})$, where $\Sigma = D \otimes _{K _{m-i}} T$. 
Clearly, exp$(\Sigma ) = p$ and deg$(\Sigma ) = p ^{m-i-1}$, so 
Brd$_{p}(K _{m-i}) \ge m - i - 1$, proving Proposition \ref{prop4.4} 
in case $i = 0$. Let finally $i > 0$. Considering inertial lifts 
over $K _{m}$ relative to $v _{i}$ of $\Sigma $ and any $L _{i} \in 
I(M _{i+1}/K _{m-i})$ with $\Sigma \otimes _{K _{m-i}} L _{i} \in 
d(L _{i})$ and $[L _{i}\colon K _{m-i}] = p ^{i}$, one obtains 
similarly that Brd$_{p}(K _{m}) \ge m - 1$.
\end{proof}
\par
\medskip
The inequalities $m - 1 \le {\rm Brd}_{p}(K) \le m$ also hold under 
the assumption that $(K, v)$ is an HDV-field, char$(K) = 0$ and 
char$(\widehat K) = p > 0$, where $\widehat K$ is an $(m - 
1)$-dimensional local field with a finite last residue field, for
some $m \ge 2$. The lower bound Brd$_{p}(K) \ge m - 1$ is obtained as
in the proof of Proposition \ref{prop4.4}, and the inequality
Brd$_{p}(K) \le m$ is implied by Proposition \ref{prop4.4} and the
injectivity of the scalar extension map Br$(K) \to {\rm
Br}(\widetilde K)$, $\widetilde K$ being the completion of $K$ with
respect to $v$ \cite{Cohn}, Theorem~1.
\par
\medskip
\section{\bf Proof of Theorem \ref{theo1.1}}

\medskip
Let $(K, v)$ be a Henselian field, $p \in \mathbb P$, $\widehat
K(p)_{\rm ab}$ the maximal abelian extension of $\widehat K$ in
$\widehat K(p)$, and $\mu _{p}(\widehat K)$, $\mu _{p}(K)$ the groups
of roots of unity of $p$-primary degrees lying in $\widehat K$ and
$K$, respectively. First, we describe index-exponent $p$-primary
$K$-pairs, assuming that $\mathcal{G}(\widehat K(p)/\widehat K)$ is a
Demushkin group and $\mu _{p}(\widehat K)$ is a nontrivial finite
group.

\medskip
\begin{lemm}
\label{lemm5.1} Let $(K, v)$ be a Henselian field containing a
primitive $p$-th root of unity, for some $p \in \mathbb P$, $p \neq
{\rm char}(\widehat K)$. Suppose that $\mathcal{G}(\widehat
K(p)/\widehat K)$ is a Demushkin pro-$p$-group, $\mu _{p}(\widehat
K)$ is a finite group of order $p ^{\nu }$, and $r _{p}(\widehat K) =
r < \infty $. Put $r ^{\prime } = r - 1$, $m ^{\prime } = {\rm
min}\{\tau (p), r ^{\prime }\}$, and for each $n \in \mathbb N$, let
$\nu _{n} = {\rm min}\{n, \nu \}$ and $\mu (p, n) = nm ^{\prime } +
\nu _{n}(m _{p} - m ^{\prime } + [(\tau (p) - m _{p})/2])$. Then $(p
^{k}, p ^{n})$, where $k, n \in \mathbb N$, is an index-exponent pair
over $K$, if and only if $n \le k \le \mu (p, n)$.
\end{lemm}

\medskip
\begin{proof}
First we prove the following assertions:
\par
\medskip
(5.1) (a) $C(\widehat K(p)/\widehat K)$ is isomorphic to the direct
sum $\mathbb Z(p ^{\infty }) ^{r'} \oplus \mathbb Z/p ^{\nu }\mathbb
Z$ and $\mathcal{G}(\widehat K(p)_{\rm ab}/\widehat K) \cong \mathbb
Z _{p} ^{r'} \oplus \mathbb Z/p ^{\nu }\mathbb Z$;
\par
(b) A cyclic extension $M$ of $\widehat K$ in $\widehat K(p)$ lies
in $I(M _{\infty }/\widehat K)$, for some $\mathbb Z _{p}$-extension
$M _{\infty }$ of $\widehat K$ in $\widehat K(p)$ if and only if
there is $M ^{\prime } \in I(\widehat K(p)/M)$, such that $M ^{\prime
}/\widehat K$ is cyclic and $[M ^{\prime }\colon M] = p ^{\nu }$;
this is the case if and only if $\mu _{p}(\widehat K) \subset N(M/\widehat
K)$.
\par
\medskip\noindent
The nontriviality of $\mu _{p}(\widehat K)$ and the Demushkin
property of $\mathcal{G}(\widehat K(p)/\widehat K)$ ensure that $r
\ge 2$, $\widehat K$ is a $p$-quasilocal nonreal field (see
\cite{Ch2}, I, Lemma~3.8). Hence, by \cite{Ch2}, I, Theorem~3.1,
Br$(\widehat K) _{p}$ is divisible, which enables one to deduce
from \cite{MS}, (11.5), and the condition on the order of $H
^{2}(\mathcal{G}(\widehat K(p)/\widehat K), \mathbb F _{p})$ that
Br$(\widehat K) _{p} \cong \mathbb Z(p ^{\infty })$. The rest of the
proof of (5.1) (a) relies on our assumption on $\mu _{p}(\widehat
K)$, which shows that $\widehat K$ contains a primitive $p ^{\nu 
}$-th root of unity $\delta $ not lying in $\widehat K ^{\ast p}$.
Consider an extension $\widehat K _{\delta }$ of $\widehat K$
obtained by adjunction of a $p$-th root of $\delta $. It is easily
verified that $\widehat K _{\delta }/\widehat K$ is a cyclic
extension of degree $p$. As $\widehat K$ is $p$-quasilocal and
Br$(\widehat K) _{p} \cong \mathbb Z(p ^{\infty })$, this means that
Br$(\widehat K _{\delta }/\widehat K)$ has order $p$. In view of
Kummer theory, cyclic $\widehat K$-algebras of degree $p$ are symbol
algebras, so the noted fact indicates that there is a cyclic degree
$p$ extension $\widehat K ^{\prime }/\widehat K$, such that the
cyclic $\widehat K$-algebra $(\widehat K ^{\prime }/\widehat K,
\sigma ', \delta )$ lies in $d(\widehat K)$ ($\sigma '$ is a
generator of $\mathcal{G}(\widehat K ^{\prime }/\widehat K)$).
Therefore, by \cite{P}, Sect. 15.1, Proposition~b, $\delta $ does
not lie in the norm group $N(\widehat K ^{\prime }/\widehat K)$.
Applying Albert's height theorem to $\widehat K ^{\prime }/\widehat
K$ (cf. \cite{FSS}, Sect. 2), one proves the nonexistence of a
cyclic extension $\widehat K _{1} ^{\prime }/\widehat K$, such that 
$[\widehat K _{1} ^{\prime }\colon \widehat K] = p ^{1+\nu }$
\par\noindent
and $\widehat K ^{\prime } \in I(\widehat K _{1} ^{\prime }/\widehat 
K)$. This result allows us to obtain from Galois theory that the 
complement $C(\widehat K(p)/\widehat K) \setminus p ^{\nu }C(\widehat 
K(p)/\widehat K)$ contains an element of order $p$. Similarly, it can 
be deduced from Kummer theory that $p ^{\nu -1}C(\widehat 
K(p)/\widehat K)$ contains all elements of $C(\widehat K(p)/\widehat 
K)$ of order $p$. Observe now that the Demushkin condition on 
$\mathcal{G}(\widehat K(p)/\widehat K)$ ensures that $C(\widehat 
K(p)/\widehat K) \cong \mathbb Z(p ^{\infty }) ^{r'} \oplus C$, for 
some cyclic $p$-group $C$ (cf. \cite{Lab}, page 106). Summing-up the 
noted properties of $C(\widehat K(p)/\widehat K)$, one concludes that 
$C \cong \mathbb Z/p ^{\nu }\mathbb Z$ and so proves (5.1) (a). As to 
(5.1) (b), it is implied by (5.1) (a) and Albert's height theorem.
\par
We continue with the proof of Lemma \ref{lemm5.1}. Statement (2.3)
(b), the isomorphism Br$(\widehat K) _{p} \cong \mathbb Z(p ^{\infty
})$, and the equality Brd$_{p}(\widehat K) = 1$ imply that $(p ^{m},
p ^{m})$, $m \in \mathbb N$, are index-exponent pairs over both
$\widehat K$ and $K$. In view of Theorem \ref{theo4.1}, this proves
Lemma \ref{lemm5.1} in the case where $\tau (p) = 1$, so we assume
that $\tau (p) \ge 2$. Suppose first that $n \in \mathbb N$ and $n
\le \nu $. Then, by Theorem \ref{theo4.1}, ind$(\Delta _{n}) \mid p
^{\mu (p,n)}$, for each $\Delta _{n} \in d(K)$ with exp$(\Delta
_{n}) = p ^{n}$. Using \cite{Mo}, Theorem~1, and the natural
bijection between $I(Y/K)$ and the set of subgroups of $v(Y)/v(K)$,
for any finite abelian tamely and totally ramified extension $Y/K$
(cf. \cite{Sch}, Ch. 3, Sect. 2), one obtains that, for each $k \in 
\mathbb N$ with $n \le k \le \mu (p, n)$, there exist an NSR-algebra 
$V _{n,k} \in d(K)$ and a totally ramified $T _{n,k} \in d(K)$, such 
that $V _{n,k} \otimes _{K} T _{n,k} \in d(K)$, exp$(V _{n,k} \otimes 
_{K} T _{n,k}) = p ^{n}$ and ind$(V _{n,k} \otimes _{K} T _{n,k}) = p 
^{k}$. These observations and the former part of (1.1) (a) prove 
Lemma \ref{lemm5.1} when $n \le \nu $. The rest of the proof is 
carried out by induction on $n \ge \nu $. The basis of the induction 
is provided by the preceding argument, which allows us to assume that 
$n > \nu $ and ind$(X) \mid p ^{\mu (p,(n-1))}$ whenever $X \in d(K)$ 
and exp$(X) \mid p ^{n-1}$. Fix an algebra $D \in d(K)$ so that 
exp$(D) = p ^{n}$ and attach to $D$ a triple $S$, $V$, $T \in d(K)$ 
as in (2.3) (a). Clearly, if exp$(V) \mid p ^{n-1}$, then exp$(V 
\otimes _{K} T) \mid p ^{n-1}$, so (4.3) and the inductive hypothesis 
imply ind$(D) \mid p ^{1+\mu (p,(n-1))} \mid p ^{\mu (p,n)}$, as 
claimed. In view of (2.4), it remains to consider the case where 
exp$(V) = p ^{n}$. Let $\Sigma $, $D _{\nu } \in d(K)$ satisfy 
$[\Sigma ] = [S \otimes _{K} V]$ and $[D _{\nu }] = p ^{\nu }[D]$ ($= 
p ^{\nu }[\Sigma ])$. Then, by (2.4) (c), exp$(\Sigma ) = p ^{n}$, 
and it follows from (4.1) and \cite{P}, Sect. 15.1, Corollary~b and 
Proposition~b, that $\Sigma /K$ is NSR. Note also that exp$(D _{\nu 
}) \mid p ^{n-\nu }$, and (2.3) (c) and \cite{P}, Sect. 15.1, 
Corollary~b, imply $D _{\nu }/K$ is NSR; in particular, $D _{\nu }$ 
contains as a maximal subfield an inertial extension $U _{\nu }$ of 
$K$. By \cite{JW}, Theorem~4.4, $U _{\nu }/K$ is abelian with 
$\mathcal{G}(U _{\nu }/K)$ of rank $u _{\nu } \le \tau (p)$. 
Moreover, it follows from (5.1), Galois theory and \cite{P}, Sect. 
15.1, Corollary~b, that $U _{\nu }$ has a $K$-isomorphic copy in $I(U 
_{\nu } ^{\prime }/K)$, for the Galois extension $U _{\nu } ^{\prime 
}$ of $K$ in $K _{\rm ur}$ with $\mathcal{G}(U _{\nu } ^{\prime }/K) 
\cong \mathbb Z _{p} ^{r'}$. Therefore, $u _{\nu } \le r ^{\prime }$, 
so \cite{JW}, Theorem~4.4, proves the following:
\par
\medskip\noindent
(5.2) ind$(D _{\nu }) \mid p ^{(n-\nu )m'}$ and $D _{\nu }$ contains
as a maximal subfield a $K$-isomorphic copy of a totally ramified
extension $\Phi _{\nu }$ of $K$ in $K(p)$.
\par
\medskip\noindent
Statement (5.2) shows that $[D _{\nu }] \in {\rm Br}(\Phi _{\nu
}/K)$, $[\Phi _{\nu }\colon K] = {\rm ind}(D _{\nu })$ and $\widehat
\Phi _{\nu } = \widehat K$. Hence, exp$(D \otimes _{K} \Phi _{\nu })
\mid p ^{\nu }$ and $r _{p}(\widehat \Phi _{\nu }) = r _{p}(\widehat
K)$, so it follows from (2.2) and Theorem \ref{theo4.1} that ind$(D
\otimes _{K} \Phi _{\nu }) \mid p ^{\nu \mu (p)}$, where $\mu (p) =
[(m _{p} + \tau (p))/2]$. As $\mu (p, n) = (n - \nu )m ^{\prime } +
\nu \mu (p)$, it is now easy to see that ind$(D) \mid p ^{\mu
(p,n)}$, as required. Suppose finally that $(k, n) \in \mathbb N
^{2}$ and $n \le k \le \mu (p, n)$. Then \cite{JW}, Example~4.3,
\cite{Mo}, Theorem~1, the above-noted properties of $U _{\nu }
^{\prime }$, and those of intermediate fields of any finite abelian 
tamely and totally ramiﬁed extension of $K$, imply the 
existence of $D _{k,n} \in d(K)$ with ind$(D _{k,n}) = p ^{k}$ and 
exp$(D _{k,n}) = p ^{n}$. Moreover, one can ensure that $D _{k,n} 
\cong N _{k,n} \otimes _{K} D _{k,n} ^{\prime }$, for some $N 
_{k,n}$, $D _{k,n} ^{\prime } \in d(K)$, such that $N _{k,n}$ is NSR 
and $D _{k,n} ^{\prime }$ is totally ramified over $K$. Lemma 
\ref{lemm5.1} is proved.
\end{proof}
\par
\medskip
Next we show that, in the setting of (1.2) (a), $C(\widehat
K(p)/\widehat K)$ possesses a divisible subgroup with infinitely many
elements of order $p$. 

\medskip
\begin{lemm}
\label{lemm5.2} Let $(E, \omega )$ be an {\rm HDV}-field with {\rm
char}$(E) = 0$, $\widehat E$ quasifinite and {\rm char}$(\widehat E) 
= p > 0$, and let $D(E(p)/E)$ be the maximal divisible subgroup of
$C(E(p)/E)$. Then:
\par
{\rm (a)} $r _{p}(E) = \infty $, provided that $\widehat E$ is
infinite;
\par
{\rm (b)} $\mu _{p}(E)$ is finite and $C(E(p)/E) \cong D(E(p)/E) 
\oplus \mathbb Z/p ^{\nu }\mathbb Z$, where $p ^{\nu }$ is the order 
of $\mu _{p}(E)$; in particular, $C(E(p)/E) = D(E(p)/E)$ if and only 
if $p ^{\nu } = 1$.
\end{lemm}
\medskip
\begin{proof}
(b): Let $\varepsilon $ be a primitive $p$-th root of unity in $E
_{\rm sep}$. It is well-known that $[E(\varepsilon )\colon E] \mid p
- 1$ (cf. \cite{L}, Ch. VIII, Sect. 3). Note also that every $E 
^{\prime } \in {\rm Fe}(E)$ is a quasilocal field with Br$(E ^{\prime 
}) \cong \mathbb Q/\mathbb Z$; hence, the scalar extension map 
\par\noindent
Br$(E) \to {\rm Br}(E ^{\prime })$ is surjective. These facts, 
combined with (1.1) (b) and \cite{P}, Sect. 15.1, Proposition~b, 
imply that if $L$ is a cyclic $p$-extension of $E$ in $E _{\rm sep}$, 
then $L(\varepsilon ) ^{\ast } = L ^{\ast }N(L(\varepsilon 
)/E(\varepsilon ))$. When $\varepsilon \notin E$, this shows that 
$\varepsilon \in N(L(\varepsilon )/E(\varepsilon ))$, which enables 
one to deduce from \cite{FSS}, Theorem~3, that $C(E(p)/E) = 
D(E(p)/E)$. Suppose now that $\mu _{p}(E) \neq \{1\}$ and denote by 
$\Gamma _{p}$ the extension of $E$ generated by the elements of $\mu 
_{p}(E _{\rm sep})$. It is known that, for any $n \in \mathbb N$, 
$\mathbb Z[X]$ contains the $p ^{n}$-th cyclotomic polynomial $\Phi 
_{p ^{n}}(X)$ (of degree $p ^{n-1}(p - 1)$), and the polynomial 
$\Phi _{p^{n}}(X + 1)$ is $p$-Eisensteinian over $\mathbb Z$. This 
implies $p ^{n-1}(p - 1)\omega _{\Gamma _{p}}(\varepsilon _{n} - 1) = 
\omega (p)$, for each $n \in \mathbb N$, $\varepsilon _{n} \in 
\Gamma _{p}$ being a primitive $p ^{n}$-th root of unity. As $\omega $ is 
discrete and $\omega (p) \neq 0$, the noted fact proves that 
$\mu _{p}(E)$ is finite. In view of \cite{Ch2}, II, Lemma~2.3, and the 
isomorphism Br$(E) _{p} \cong \mathbb Z(p ^{\infty })$, the obtained 
result yields $C(E(p)/E) \cong D(E(p)/E) \oplus \mathbb Z/p ^{\nu
}\mathbb Z$, as claimed by Lemma \ref{lemm5.2} (b).
\par
(a): Assume that $\widehat E$ is infinite, fix a uniformizer $\pi
\in E$ and elements $a _{n} \in E$, $n \in \mathbb N$, so that
$\omega (a _{n}) = 0$ and the residue classes $\hat a _{n}$, $n
\in \mathbb N$, be linearly independent over the prime subfield
$\mathbb F _{p}$ of $\widehat E$. It is easily verified that the
cosets $(1 + a _{n}\pi )E ^{\ast p}$, $n \in \mathbb N$, are
linearly independent over $\mathbb F _{p}$. This means that $E
^{\ast }/E ^{\ast p}$ is an infinite group. At the same time, by
local class field theory, if $L _{1}, \dots , L _{n}$ are cyclic
extensions of $E$ in $E(p)$ of degree $p$, and $L = L _{1} \dots L
_{n}$, then $E ^{\ast p} \le N(L/E) \le E ^{\ast }$ and the index of
$N(L/E)$ in $E ^{\ast }$ is equal to $[L\colon E]$. Finally, the
quasilocality of $E$ shows that if $a \in E ^{\ast } \setminus E
^{\ast p}$, $D \in d(E)$ and ind$(D) = p$, then there is a cyclic
degree $p$ extension $Y$ of $E$ in $E(p)$, such that $D \cong (Y/E,
\tau , a)$, for some generator $\tau $ of $\mathcal{G}(Y/E)$ (cf.
\cite{P}, Sect. 15.5, and \cite{Ch2}, I, Corollary~8.5). Hence, by
\cite{P}, Sect. 15.1, Proposition~b, $a \notin N(Y/E)$, which means
that $E ^{\ast p}$ equals the intersection of the norm groups of
cyclic extensions of $E$ of degree $p$. Now it is clear that $r 
_{p}(E) = \infty $, so Lemma \ref{lemm5.2} is proved.
\end{proof}
\par
\medskip
{\it We are now in a position to prove (1.2) (a)}. The fulfillment of
the conditions of Lemma \ref{lemm5.2} ensures that $D(E(p)/E)$
contains infinitely many elements of order $p$. Hence, by Galois
theory and the divisibility of $D(E(p)/E)$, every finite abelian
$p$-group $G$ is isomorphic to a subgroup of $D(E(p)/E)$. Assuming
now that $E$ is isomorphic to $\widehat K$, for some Henselian field
$(K, v)$, and using \cite{TW}, Theorem~A.24, one obtains further that
$K$ possesses a Galois extension $U _{G}$ in $K _{\rm ur}$ with
$\mathcal{G}(U _{G}/K) \cong G$. When the rank of $G$ is at most
$\tau (p)$, one deduces from \cite{Mo}, Theorem~1 (or \cite{JW},
Example~4.3), that there is an NSR-algebra $D _{G} \in d(K)$
possessing a maximal subfield $K$-isomorphic to $U _{G}$. Thus it
becomes clear that there exist $D _{k,n} \in d(K)\colon k, n \in
\mathbb N, n \le k \le \tau (p)n$, such that $D _{k,n}/K$ is NSR,
ind$(D _{k,n}) = p ^{k}$ and exp$(D _{k,n}) = p ^{n}$. The obtained
result proves (1.2) (a), since Theorem \ref{theo4.1} and the equality
$r _{p}(E) = r _{p}(\widehat K) = \infty $ yield Brd$_{p}(K) = \tau
(p)$.
\par
\medskip
{\it Our objective now is to prove (1.2) (b), (c) and (d)}. Suppose 
that $(K, v)$ is Henselian, such that $v(K) \neq pv(K)$, Brd$_{p}(K) 
< \infty $, and $\widehat K$ has a Henselian discrete valuation 
$\omega $ whose residue field $\tilde k$ is quasifinite with 
char$(\tilde k) \neq p$. Then $\widehat K$ is quasilocal and
Brd$_{p}(K)$ is determined by Theorem \ref{theo4.1} (a). Also, the 
conditions on $\omega $ ensure that $\widehat K ^{\ast }/\widehat K
^{\ast p} \cong \tilde k ^{\ast }/\tilde k ^{\ast p} \times \omega 
(\widehat K)/p\omega (\widehat K)$. This allows to prove those of 
the following statements, for which we assume that $\mu 
_{p}(\widehat K) \neq \{1\}$:
\par
\medskip
(5.3) (a) $r _{p}(\widehat K) \le 2$ and $r _{p}(\widehat K) = 2
\leftrightarrow \mu _{p}(\widehat K) \neq \{1\}$ (cf. \cite{FV}, Ch.
2, (3.5));
\par
(b) If $\mu _{p}(\widehat K) = \{1\}$, then finite extensions of
$\widehat K$ in $\widehat K(p)$ are inertial relative to $\omega $,
and $\mathcal{G}(\widehat K(p)/\widehat K) \cong \mathcal{G}(\tilde
k(p)/\tilde k) \cong \mathbb Z _{p}$ (see \cite{Wh}, Theorem~2, and
\cite{Ch1}, Lemma~1.1);
\par
(c) $\mathcal{G}(\widehat K(p)/\widehat K)$ is a Demushkin group
when $\mu _{p}(\tilde k) \neq \{1\}$ (cf. \cite{Wa}, Lemma~7);
\par
(d) $\mathcal{G}(\widehat K _{\rm ab}(p)/\widehat K) \cong \mathbb
Z _{p} \oplus \mathbb Z/p ^{\nu }\mathbb Z$, provided that $\mu 
_{p}(\tilde k)$ is of finite order $p ^{\nu }$; 
$\mathcal{G}(\widehat K _{\rm ab}(p)/\widehat K) \cong \mathbb Z 
_{p}^{2}$, if $\mu _{p}(\tilde k)$ is infinite (apply (5.1) (a) in 
the former case, and use Kummer theory in the latter one). 
\par
\medskip
The inequality $p \neq {\rm char}(\tilde k)$ and the quasilocality of
$\widehat K$ show that Brd$_{p}(K)$ can be determined by applying 
Theorem \ref{theo4.1}. In view of (5.3) (a), (b) and the 
divisibility of Br$(\widehat K) _{p}$, this proves (1.2) (b) and 
(c). The former part of (1.2) (d) follows from (5.3) (c), (d) and 
Lemma \ref{lemm5.1}; in this case, $\mu (p, n)$ is equal to $n + 
{\rm min}\{n, \nu \}[\tau (p)/2]$, for each $n \in \mathbb N$. For 
the proof of the latter one, we use the concluding part of (5.3) 
(d), which implies every finite abelian $p$-group $G$ of rank $\le 
2$ is isomorphic to $\mathcal{G}(U _{G}/K)$, for some Galois 
extension $U _{G}$ of $K$ in $K _{\rm ur}$. This gives us the 
possibility to complete the proof of (1.2) (d), arguing along the 
lines drawn at the end of the proof of (1.2) (a).
\par
\medskip
{\it We prove Theorem \ref{theo1.1}}. The field $\widehat K$ is
quasilocal, and is complete relative to a discrete valuation $\omega
$ with a finite residue field $\tilde k$. This implies $\omega $ is
Henselian, $\mu _{p}(\widehat K)$ is finite, Br$(\widehat K) \cong
\mathbb Q/\mathbb Z$, and in case $p \neq {\rm char}(\tilde k)$,
$\varepsilon _{p} \in \widehat K$ if and only if $p$ divides the
order of $\tilde k ^{\ast }$. When $\varepsilon _{p} \notin \widehat 
K$, $C(\widehat K(p)/\widehat K)$ is divisible, by the following 
results (which are contained in (5.3) (b) and \cite{S2}, Theorem~3, 
respectively):
\par
\medskip
(5.4) (a) $\mathcal{G}(\widehat K(p)/\widehat K) \cong \mathbb Z
_{p}$, provided that $p \neq {\rm char}(\tilde k)$;
\par
(b) If char$(\widehat K) = 0$ and char$(\tilde k) = p$, then
$\mathcal{G}(\widehat K(p)/\widehat K)$ is a free pro-$p$-group, and
$\mathcal{G}(\widehat K(p)_{\rm ab}/\widehat K) \cong \mathbb Z _{p}
^{r}$, where $r = r _{p}(\widehat K)$; in addition, $\widehat K$ is
a finite extension of the field $\mathbb Q _{p}$ of $p$-adic numbers
and $r = [\widehat K\colon \mathbb Q _{p}] + 1$.
\par
\medskip\noindent
Note also that, by Theorem \ref{theo4.1}, Brd$_{p}(K) = m _{p}$, and
by (5.4) and \cite{Mo}, Theorem~1, each pair of $p$-primary integers
admissible by Theorem \ref{theo1.1} is an index-exponent pair of a
suitably chosen NSR-algebra over $K$.
\par
Consider finally the case where $\varepsilon _{p} \in \widehat K$.
Then Theorem \ref{theo4.1} yields Brd$_{p}(K) = \mu (p, 1)$, and
Lemma \ref{lemm5.1} implies $(1, 1)$ and $(p ^{k}, p ^{n})\colon k,
n \in \mathbb N, n \le k \le \mu (p, n)$, are all index-exponent
$p$-primary $K$-pairs. This completes our proof.

\medskip
\begin{rema}
\label{rema5.3}
Theorem \ref{theo1.1} retains validity, if $\widehat K \in {\rm
Fe}(\mathbb Q _{\pi } ^{\prime })$, for some $\pi $-adically closed
field $\mathbb Q _{\pi } ^{\prime }$ (in the sense of \cite{PR}).
This is fulfilled, if char$(\widehat K) = 0$ and $\widehat K$ has a
Henselian discrete valuation $\omega $ with a finite residue field
$\tilde k$ of characteristic $\pi $. Also, (5.4) hold, if $\mu
_{p}(\widehat K) = \{1\}$ (in case (b), with $\mathbb Q _{p}
^{\prime }$ instead of $\mathbb Q _{p}$). When $\mu _{p}(\widehat K)
\neq \{1\}$ and $r = r _{p}(\widehat K)$, we have: $r = 2$, provided
$p \neq \pi $; $r = [\widehat K\colon \mathbb Q _{p} ^{\prime }] +
2$, if $p = \pi $ (see (5.3), \cite{Wa}, Lemma~7, and \cite{Lab}, 
Sect. 5, for the case of $\mathbb Q _{p} ^{\prime } = \mathbb Q 
_{p}$).
\end{rema}

\medskip
\begin{coro}
\label{coro5.4} Let $(K, v)$ be a Henselian field, such that $\tau 
(p) < \infty $, for some $p \in \mathbb P$, $p \neq {\rm 
char}(\widehat K)$. Also, let $\widehat K$ have a Henselian discrete 
valuation $\omega $ with a quasifinite residue field $\tilde k$. Then 
{\rm abrd}$_{p}(K) = 1 + [\tau (p)/2]$, if $p \neq {\rm 
char}(\tilde k)$; {\rm abrd}$_{p}(K) = {\rm max}\{1, \tau (p)\}$, if
{\rm char}$(\widehat K) = 0$ and {\rm char}$(\tilde k) = p$.
\end{coro}

\medskip
\begin{proof}
In view of (1.1) (b) and (1.2), one may consider only the case
where $\mu _{p}(\widehat K) \neq  \{1\}$, char$(\widehat K) = 0$,
$\tilde k$ is finite and char$(\tilde k) = p$. Then our conclusion
follows from Remark \ref{rema5.3} and the fact that $[\widehat
K(p)\colon \widehat K] = \infty $.
\end{proof}

\medskip
\emph{Conclusion.} Assume that $(K, v)$ is Henselian with $\widehat 
K$ possessing a Henselian discrete valuation $\omega $ whose residue 
field is quasifinite. Summing-up (1.1), (2.3) (b) and Corollary 
\ref{coro2.2}, observing that Br$(\widehat K) \cong \mathbb Q/\mathbb 
Z$ and Brd$_{p}(\widehat K) = 1$, $p \in \mathbb P$, and using 
results of this paper, one describes index-exponent $K$-pairs 
prime-to char$(\widehat K)$. The non-divisibility restriction is 
superfluous, if char$(K) > 0$, $(K, v)$ is maximally complete and 
$\widehat K$ satisfies the conditions of Corollary \ref{coro3.6}.

\vskip0.2truecm
\emph{Acknowledgement.} The present research was partially supported
by Project No. RD-08-118/04.02.2019 of Shumen University, Bulgaria.

\end{document}